\documentclass[onefignum,onetabnum]{siamart190516}

\usepackage{lipsum}
\usepackage{amsfonts}
\usepackage{graphicx}
\usepackage{epstopdf}
\graphicspath{{img/}}
\ifpdf
  \DeclareGraphicsExtensions{.eps,.pdf,.png,.jpg}
\else
  \DeclareGraphicsExtensions{.eps}
\fi

\usepackage{booktabs,caption}
\usepackage[flushleft]{threeparttable}

\usepackage{amsmath}
\usepackage{amssymb}
\usepackage{url}
\usepackage{listings}
\usepackage{color}
\usepackage{algorithm}
\usepackage{algorithmicx}
\usepackage{algpseudocode}
\usepackage{verbatim}
\usepackage{tikz}
\usepackage{tikz-qtree}
\usepackage{wasysym}
\usepackage{mathtools}
\usepackage{subfig}

\newcommand{\bigo}{\mathcal{O}}
\newcommand{\isep}{\mathrel{{.}\,{.}}\nobreak}

\algnewcommand\algorithmicinput{\textbf{Input:}}
\algnewcommand\INPUT{\item[\algorithmicinput]}
\algnewcommand\algorithmicoutput{\textbf{Output:}}
\algnewcommand\OUTPUT{\item[\algorithmicoutput]}
\algnewcommand{\IIf}[1]{\State\algorithmicif\ #1\ \algorithmicthen}
\algnewcommand{\EndIIf}{\unskip\ \algorithmicend\ \algorithmicif}

\definecolor{darkcyan}{rgb}{0.0, 0.55, 0.55}
\definecolor{codegreen}{rgb}{0,0.6,0}
\definecolor{codegray}{rgb}{0.5,0.5,0.5}
\definecolor{codepurple}{rgb}{0.58,0,0.82}
\definecolor{backcolour}{rgb}{0.95,0.95,0.92}

\newsiamremark{remark}{Remark}
\newsiamremark{hypothesis}{Hypothesis}
\crefname{hypothesis}{Hypothesis}{Hypotheses}
\newsiamthm{claim}{Claim}

\headers{Algorithms minimizing the max and sum of vectors of functions}{H. Khaleghzadeh, R. R. Manumachu, and A. Lastovetsky}

\title{Novel bi-objective optimization algorithms minimizing the max and sum of vectors of  functions\thanks{\funding{This publication has emanated from research conducted with the financial support of Science Foundation Ireland (SFI) under Grant Number 14/IA/2474.}}}

\author{Hamidreza Khaleghzadeh\thanks{School of Computer Science, University College Dublin, Belfield, Dublin 4, Ireland 
  (\email{hamidreza.khaleghzadeh@ucd.ie}, \email{ravi.manumachu@ucd.ie}, \email{alexey.lastovetsky@ucd.ie})}
\and Ravi~Reddy~Manumachu\footnotemark[2]
\and Alexey Lastovetsky\footnotemark[2]}

\usepackage{amsopn}

\ifpdf
\hypersetup{
  pdftitle={Algorithms minimizing the max and sum of vector of functions},
  pdfauthor={H. Khaleghzadeh, R. R. Manumachu, and A. Lastovetsky}
}
\fi

\usepackage[utf8]{inputenc}
\usepackage[T1]{fontenc}

\begin{document}

\maketitle

\begin{abstract}
We study a bi-objective optimization problem, which for a given positive real number $n$ aims to find a vector $X = \{x_0,\cdots,x_{k-1}\} \in  \mathbb{R}^{k}_{\ge 0}$ such that  $\sum_{i=0}^{k-1} x_i = n$, minimizing the maximum of $k$  functions of objective type one, $\max_{i=0}^{k-1} f_i(x_i)$, and the sum of $k$  functions of objective type two, $\sum_{i=0}^{k-1} g_i(x_i)$. This problem arises in the optimization of applications for performance and energy on high performance computing platforms. We first propose an algorithm solving the problem for the case where all the functions of objective type one are continuous and strictly increasing, and all the functions of objective type two are linear increasing. We then propose an algorithm solving a version of the problem where $n$ is a positive integer and all the functions are discrete and represented by finite sets with no assumption on their shapes. Both algorithms are of polynomial complexity.
\end{abstract}



\section{Introduction}

Bi-objective optimization problems where one objective function is \textit{max} and the other is \textit{sum} are common in the category of multiple objective minimum spanning tree problems that have important applications in the field of network design and optimization (transportation and communication networks, for example) \cite{Hansen1980}, \cite{Berman1990}, \cite{minoux1989solving}, \cite{Punnen1994}, \cite{Punnen1996}, \cite{sergienko1987finding}, \cite{melamed1996computational}, \cite{melamed1998numerical}, \cite{Leizer2009}, \cite{Leizer2010}, \cite{Bornstein2012}, \cite{Pinto2019}, flowshop group scheduling \cite{karimi2010bi},\cite{Torkashvand2017},  after disaster blood supply chain management \cite{heydari2018dynamic}, collaborative production planning \cite{salamati2018trade}, and performance and energy optimization of high computing systems and applications \cite{tarplee2016energy}, \cite{aba2017approximation}, \cite{manumachu2018bi}, \cite{manumachu2018bicpe}, \cite{HamidTPDS2020}.

In this work, we introduce two mathematical problems motivated by the problem of bi-objective optimization of scientific applications on modern heterogeneous high-performance computing (HPC) platforms for performance and energy.

To motivate the first problem, consider the bi-objective optimization of a popular and highly optimized matrix multiplication application on a hybrid heterogeneous computing platform for performance and energy. The optimization goal is to find performance-energy optimal application configurations (workload distributions), minimizing the execution time (\emph{min-max}) and the total energy consumption (\emph{min-sum}) of computations during the parallel execution of the application.

The application computes the matrix product, $C = \alpha \times A \times B + \beta \times C$, where $A$, $B$, and $C$ are matrices of size $M \times N$, $N \times N$, and $M \times N$, and $\alpha$ and $\beta$ are constant floating-point numbers. The platform consists of five heterogeneous processors: Intel Haswell E5-2670V3 multi-core CPU (CPU\_1), Intel Xeon Gold 6152 multi-core CPU (CPU\_2), NVIDIA K40c GPU (GPU\_1), NVIDIA P100 PCIe GPU (GPU\_2), and Intel Xeon Phi 3120P (XeonPhi\_1). Further details of our computing platform are given in the Appendix (Section \ref{app:platform}).

\begin{figure}[!htbp]
	\centering
	\centering
	\subfloat[]{{ \includegraphics[height=2.3in,width=4in]{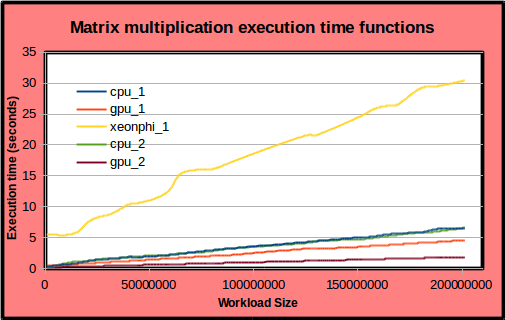} } \label{fig:dgemm_time}}
	\quad
	\subfloat[]{{ \includegraphics[height=2.3in,width=4in]{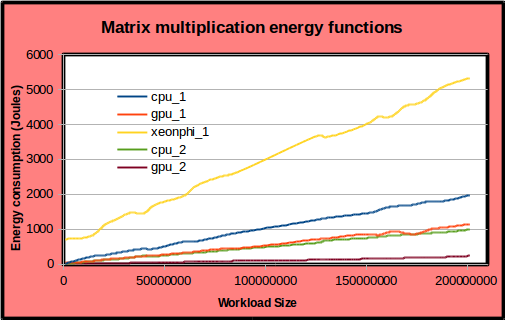} } \label{fig:dgemm_eng}}
	\caption{Execution time and energy functions of matrix multiplication application.}
	\label{fig:dgemm_functions}
\end{figure}

The figures \ref{fig:dgemm_time} and \ref{fig:dgemm_eng} show the execution time functions $\{f_0(x),\dots,f_4(x)\}$ and the energy functions $\{g_0(x),\dots,g_4(x)\}$ of the processors against the workload size ($x$). The energy consumption during an application execution is obtained using system-level physical power measurements using power meters, which is considered the most accurate method of energy measurement \cite{fahad2019comparative}. The execution time function shapes are continuous and strictly increasing. The energy function shapes can be approximated accurately by linear increasing functions. The optimization goal is to find workload distributions of the workload size $n$ ($\{x_0,\dots,x_4\}, \sum_{i=0}^4 x_i = n$) minimizing the execution time ($\max_{i=0}^4 f_i(x_i)$) and the total energy consumption ($\sum_{i=0}^4 g_i(x_i)$) during the parallel execution of the application.

To the best of our knowledge, there is no research work tackling bi-objective optimization problems aiming to minimize the max of $k$-dimensional vector of functions of objective type one and the sum of $k$-dimensional vector of functions of objective type two subject to linear constraints. $k$ is equal to 5 in the example illustrated above. State-of-the-art bi-objective optimization methods consider objective functions and constraints that are linear functions of the decision variables (except \cite{Pinto2019} who consider a non-linear function for the max objective). The objective functions are typically \textit{max} of a function of objective type one and \textit{sum} of a function of objective type two. Therefore, the state-of-the-art methods do not consider vectors of functions. 

In this work, we formulate the mathematical problem, which for a given positive real number $n$ aims to find a vector $X = \{x_0,\cdots,x_{k-1}\} \in  \mathbb{R}^{k}_{\ge 0}$ such that  $\sum_{i=0}^{k-1} x_i = n$, minimizing the max of $k$-dimensional vector of  functions of objective type one and the sum of $k$-dimensional vector of functions of objective type two. We propose an algorithm solving the case where all the functions of objective type one are continuous and strictly increasing, and all the functions of objective type two are linear increasing. The algorithm exhibits polynomial complexity.

\begin{figure}[!htbp]
	\centering
	\centering
	\subfloat[]{{ \includegraphics[height=2.3in,width=4in]{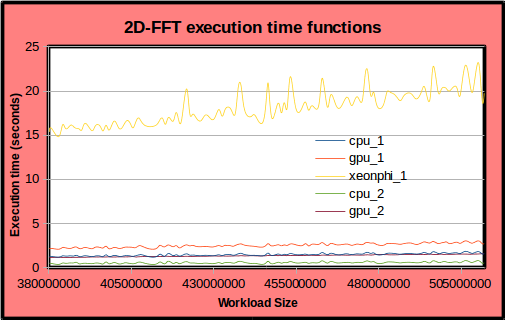} } \label{fig:fft_time}}
	\quad
	\subfloat[]{{ \includegraphics[height=2.3in,width=4in]{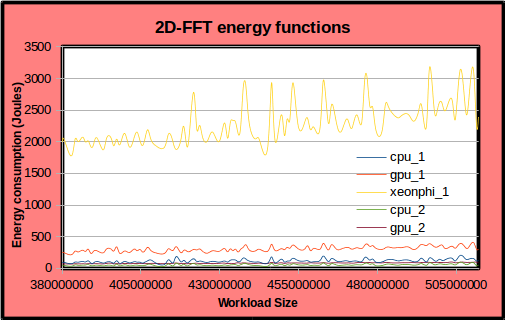} } \label{fig:fft_eng}}
	\caption{Execution time and energy functions of 2D-FFT application.}
	\label{fig:fft_functions}
\end{figure}

To motivate the second problem, consider the bi-objective optimization of a popular and highly optimized 2D fast Fourier transform application (2D-FFT) on the same platform for performance and energy. The application computes 2D discrete Fourier transform of a complex signal matrix of size $M \times N$. The figures \ref{fig:fft_time} and \ref{fig:fft_eng} show the execution time functions $\{f_0(x),\dots,f_4(x)\}$ and the energy functions $\{g_0(x),\dots,g_4(x)\}$ of the processors against the workload size ($x$). While the execution time function shapes of the processors (CPU\_1, GPU\_1, CPU\_2, GPU\_2) can be approximated accurately by continuous and strictly increasing functions, the shape for the processor XeonPhi\_1 is not amenable to such approximation. Similarly, the shape for the processor XeonPhi\_1 is not amenable to linear approximation, whereas the energy function shapes of the processors (CPU\_1, CPU\_2, GPU\_1, GPU\_2) can be approximated accurately by linear increasing functions.

We formulate the mathematical problem dealing with this scenario. Given a positive integer, $n$, the problem aims to find a vector, $X = \{x_0,\cdots,x_{k-1}\} \in \mathbb{Z}^{k}_{\ge 0}$, such that  $\sum_{i=0}^{k-1} x_i = n$, minimizing the max of $k$-dimensional vector of discrete functions of objective type one and the sum of $k$-dimensional vector of discrete functions of objective type two. The functions are represented by finite sets of cardinality $m$. We propose a branch-and-bound algorithm that exhibits polynomial complexity to solve the problem when the functions are of an arbitrary shape.

The main original contributions of this work are:

\begin{itemize}
	\item Mathematical formulation of the bi-objective optimization problem which for a given positive real number $n$ aims to find a vector, $X = \{x_0,\cdots,x_{k-1}\} \in \mathbb{R}^{k}_{\ge 0}$, such that  $\sum_{i=0}^{k-1} x_i = n$, minimizing the maximum of $k$ functions of objective type one and the sum of $k$ functions of objective type two.
	\item An exact algorithm of polynomial complexity solving the bi-objective optimization problem when all the functions of objective type one are continuous and strictly increasing, and all the functions of objective type two are linear increasing.
	\item Mathematical formulation of the bi-objective optimization problem which for a given positive integer, $n$, aims to find a vector, $X = \{x_0,\cdots,x_{k-1}\} \in \mathbb{Z}^{k}_{\ge 0}$, such that $\sum_{i=0}^{k-1} x_i = n$, minimizing the maximum of $k$ discrete functions of objective type one and the sum of $k$ discrete functions of objective type two. The functions are represented by finite sets of cardinality $m$. 
    \item An exact algorithm of polynomial complexity solving the bi-objective optimization problem when the functions are of an arbitrary shape.
\end{itemize}

The rest of the paper is organized as follows. We discuss the related work in section \ref{sec:related-work}. The formulation of the bi-objective optimization problem is presented in section \ref{sec:bop_formulation}. In section \ref{sec:lbopa}, we propose an algorithm solving the bi-objective optimization problem when all the functions of objective type one are continuous and strictly increasing, and all the functions of objective type two are linear increasing. Section \ref{sec:gbopa} presents an algorithm solving the problem where all the functions are discrete functions of an arbitrary shape. Finally, we conclude the paper in section \ref{sec:conclusion}. 

\section{Related work} \label{sec:related-work}

Our survey of related work focuses exclusively on research works that consider optimization problems where the number of objectives is at least two and where both the types of objectives (max and sum) are present. We also cover optimization problems that minimize the sum of the two types of objectives. For each research work, we specify whether the solution method is an exact method or a metaheuristic. Since the focus of this work is bi-objective optimization, we define the problem first.

In mathematical terms, a bi-objective optimization problem can be formulated as \cite{Miettinen1999},\cite{Talbi2009}:
\begin{alignat*}{3}
	& minimize \quad \{T(\vec{x}),E(\vec{x})\} & \\
	& \text{Subject to} \quad \vec{x} \in \mathcal{X}
\end{alignat*}

There are two objective functions, $T:\mathbb{R}^k \rightarrow \mathbb{R}$ and $E:\mathbb{R}^k \rightarrow \mathbb{R}$ denoted by a vector, $\mathcal{F}(\vec{x})=(T(\vec{x}),E(\vec{x}))^T$. The decision variable vectors, $\vec{x} = (x_1,...,x_k)^T$, belong to the feasible set, $\mathcal{X} \subset \mathbb{R}^k$. The feasible objective set, which is the image of $\mathcal{X}$, is given by $\mathcal{F}(\mathcal{X}) \subset \mathbb{R}^2$. The goal of the problem is to minimize both the objective functions simultaneously. The solutions of the problem are called Pareto-optimal solutions. The set of Pareto-optimal solutions is called the Pareto front.

\begin{definition}
A decision variable vector $\vec{x}^* \in \mathcal{X}$ is \textit{Pareto optimal} if there does not exist another decision variable vector $\vec{x} \in \mathcal{X}$ such that $T(\vec{x}) \leq T(\vec{x}^*), E(\vec{x}) \leq E(\vec{x}^*)$ and either $T(\vec{x}) < T(\vec{x}^*)$ or $E(\vec{x}) < E(\vec{x}^*)$ or both.
\end{definition}

Any point in the feasible objective set, $\mathcal{F}(\mathcal{X})$, that is not on the Pareto front is a bad solution. Along the Pareto front, moving from one Pareto-optimal solution to the other will require a trade-off between the two objectives.

There are several classifications for methods solving bi-objective optimization problems \cite{Miettinen1999},\cite{Talbi2009}. Since the set of Pareto optimal solutions is partially ordered, one classification is based on the involvement of the decision-maker in the solution method to select specific solutions. There are four categories in this classification, \textit{No preference}, \textit{A priori}, \textit{A posteriori}, \textit{Interactive}, which are described in the Appendix (Section \ref{app:solmethods}). The algorithms solving bi-objective optimization problems can be divided into two major categories, \textit{exact methods} and \textit{metaheuristics}. While branch-and-bound (B\&B) is the dominant technique in the first category, genetic algorithm (GA) is popular in the second category.

\subsection{Mathematical Multi-Objective Optimization}

Hansen \cite{Hansen1980} study bi-objective optimization problems related to path selection (shortest path, minimum spanning tree) for directed graphs. Berman et al. \cite{Berman1990} improve the solutions proposed by \cite{Hansen1980}. 

Minoux \cite{minoux1989solving}, Punnen \cite{Punnen1994}, Punnen and Nair \cite{Punnen1996} study a class of combinatorial optimization problems in which the objective function minimized is an algebraic sum of a bottleneck cost function (Min-Max) and a linear cost function (Min-Sum). The authors present a solution method that employs weighting method (where the weighting coefficients are 1) to find a good compromise solution between the two objectives. Sergienko and Perepelitsa \cite{sergienko1987finding}, Melamed and Sigal \cite{melamed1996computational} present a solution method specifically designed to solve optimization problems containing one max objective and one sum objective. Melamed and Sigal \cite{melamed1998numerical} propose a solution method specifically designed to solve optimization problems containing two max objectives and one sum objective. Ruzika and Hamacher \cite{Ruzika2009} present a survey on minimum spanning tree problems with two or more objective functions that can be \textit{max} or \textit{sum}. The proposed solution methods \cite{sergienko1987finding}, \cite{melamed1996computational}, \cite{melamed1998numerical} are exact and return the optimal set of solutions.

Leizer et al. \cite{Leizer2009}, \cite{Leizer2010} propose a solution method solving a tri-objective path selection problem involving two bottleneck functions and a cost function. Bornstein et al. \cite{Bornstein2012} consider optimal spanning trees and optimal paths problems involving one cost function and one or more bottleneck functions. Their algorithm determines the Pareto front in polynomial time if the cost function is solvable in polynomial time. Pinto et al. \cite{Pinto2019} present a bi-objective network flow routing problem regarding network load balancing and flow path length. They propose an exact and polynomial approach. In their formulation, the cost function and the constraints are linear and the bottleneck function is non-linear.

\subsection{Domain-specific Bi-Objective Optimization}

Multi-objective flow shop scheduling deals with scheduling jobs on machines to optimize two or more objectives. The flow shop scheduling problem is NP-Complete, and therefore metaheuristics are heavily employed in solution methods for the problem. Karimi et al. \cite{karimi2010bi} propose a genetic algorithm to solve a flow shop scheduling problem where objective functions are makespan (max) and total weighted tardiness (sum). Sun et al. \cite{Sun2011} present a taxonomy of research in multi-objective flow shop scheduling. Torkashvand et al. \cite{Torkashvand2017} study multi-objective flow shop scheduling problems with interfering jobs where there are two sets of jobs, each with its objective. Some jobs are scheduled to minimize makespan while others are to minimize total tardiness.

Heydari et al. \cite{heydari2018dynamic} propose a bi-objective model for after disaster blood supply chain management. The problem is formulated as a mixed-integer linear programming model and considers two objectives, the amount of blood shortage in demand points (max) and total cost, which is the summation of costs for blood transportation, process, and holding (sum). 

Salamati-Hormozi et al. \cite{salamati2018trade} survey bi-objective optimization problems involving minimization of the total cost (sum) and minimization of the
maximal production utilization to achieve fair allocations of production loads (max). The authors propose metaheuristic algorithms and the $\epsilon$-constraint method to determine the Pareto front. 

\subsubsection{Bi-Objective Optimization on High Performance Computing Platforms}

There are two principal categories of solution methods for optimizing applications on high performance computing (HPC) platforms for performance and energy. The first category of system-level solution methods aim to optimize the performance and energy of the executing environment of the applications. The dominant decision variable in this category is Dynamic Voltage and Frequency Scaling (DVFS). DVFS reduces the dynamic power consumed by a processor by throttling its clock frequency. The methods proposed in \cite{Yu2015},\cite{gholkar2016power},\cite{Rountree2017} optimize for performance under a energy budget or optimize for energy under an execution time constraint. The methods proposed in \cite{Kessaci2013},\cite{Durillo2014},\cite{Kolodziej2015} solve bi-objective optimization for performance and energy with no time constraint or energy budget.

The second category of application-level solution methods \cite{Lang2014},\cite{chakrabarti2017pareto},\cite{LastovetskyReddy2017},\cite{manumachu2018bi},\cite{manumachu2018bicpe},\cite{HamidTPDS2020} use application-level decision variables and models. The most popular decision variables include the loop tile size, workload distribution, number of processors, and number of threads.

Reddy et al. \cite{manumachu2018bi}, \cite{manumachu2018bicpe} study bi-objective optimization of data-parallel applications for performance and energy on homogeneous clusters multicore CPUs employing only one decision variable, the workload distribution. They propose an efficient solution method. The method accepts as input the number of available processors, the discrete function of the processor's energy consumption against the workload size, the discrete function of the processor's performance against the workload size. It outputs a Pareto-optimal set of workload distributions. Khaleghzadeh et al. \cite{HamidTPDS2020} propose exact solution methods solving bi-objective optimization problem for hybrid data-parallel applications on heterogeneous computing platforms for performance and energy.

Tarplee et al. \cite{tarplee2016energy} consider optimizing two conflicting objectives, the make-span and total energy consumption of all nodes in a HPC platform. They employ linear programming and divisible load theory to compute tight lower bounds on the make-span and energy of all tasks on a given platform. Using this formulation, they then generate a set of Pareto front solutions. The decision variable is task mapping. Aba et al. \cite{aba2017approximation} present an approximation algorithm to minimize both make-span and the total energy consumption in parallel applications running on a heterogeneous resources system. The decision variable is task scheduling. Their algorithm ignores all solutions where energy consumption exceeds a given constraint and returns the solution with minimum execution time.

\subsection{Summary}

State-of-the-art bi-objective optimization methods consider objective functions and constraints that are linear functions of the decision variables (except \cite{Pinto2019} who consider a non-linear function for the max objective). The objective functions are typically max of a function of objective type one and sum of a function of objective type two. To the best of our knowledge, there is no research work tackling bi-objective optimization problems aiming to minimize the max of $k$-dimensional vector of functions of objective type one and the sum of $k$-dimensional vector of functions of objective type two subject to linear constraints. This problem arises in optimizing applications for performance and energy consumption on modern heterogeneous HPC platforms \cite{HamidTPDS2020}. We address the gap in this work.

\section{Formulation of the Bi-objective Optimization Problem} \label{sec:bop_formulation}

Given a positive real number $n \in \mathbb{R}_{> 0}$  and two sets of  $k$  functions each, $F = \{f_0, f_1, \cdots, f_{k-1}\}$ and  $G = \{g_0, g_1, \cdots, g_{k-1}\}$, where $f_i, g_i \colon \mathbb{R}_{\ge 0} \to \mathbb{R}_{\ge 0}, i \in \{0,\cdots,k-1\}$, the problem is to find a vector $X = \{x_0,\cdots,x_{k-1}\} \in  \mathbb{R}^{k}_{\ge 0}$ such that $\sum_{i=0}^{k-1} x_i = n$, minimizing the objective functions $T(X) = \max_{i=0}^{k-1}~f_i(x_i)$ and $ E(X) =  \sum_{i=0}^{k-1} g_i(x_i)$. We use $T \times E$ to denote the objective space of this problem, $\mathbb{R}_{\ge 0} \times \mathbb{R}_{\ge 0}$.

Thus, the problem can be formulated as follows:

\textbf{BOPGVEC($n, k, F, G$)}:
\begin{equation} \label{eq:bopfg}
\begin{split}
T(X) &= \max_{i=0}^{k-1} f_i(x_i) \\
E(X) &= \sum_{i=0}^{k-1} g_i(x_i) \\
\underset{X}{\text{minimize}} & \quad \{T(X), E(X)\} \\
\quad \text{s.t.} & \quad x_0 + x_1 + \cdots + x_{k-1} = n
\end{split}
\end{equation}

We aim to solve BOPGVEC by finding both the Pareto front containing the optimal objective vectors in the objective space $T \times E$ and the decision vector for a point in the Pareto front. Thus, our solution finds a set of triplets $\Psi = \{(T(X), E(X), X)\}$ such that $X$ is a Pareto-optimal  decision vector, and the projection of $\Psi$ onto the objective space $T \times E$, $\Psi\downarrow_{T \times E}$, is the Pareto front.

\section{Bi-objective Optimization Problem for Max of Increasing Continuous Functions and Sum of Linear Increasing Functions} \label{sec:lbopa}

In this section, we solve BOPGVEC for the case where all functions in the set $F$ are continuous and strictly increasing, and all functions in the set $G$ are linear increasing, that is, $G = \{g_0, \cdots, g_{k-1}\}, g_i(x) = b_i \times x,  b_i \in \mathbb{R}_{> 0}, i = 0,\dots,k-1$. Without loss of generality, we assume that the functions in $G$ are sorted in the decreasing order of coefficients, $b_0 \ge b_1 \ge \dots \ge b_{k-1}$.

Our solution consists of two algorithms, Algorithm \ref{alg:lbopa} and Algorithm \ref{alg:partition}. The first one, which we call LBOPA, constructs the Pareto front of the optimal solutions in the objective space  $\Psi\downarrow_{T \times E}$. The second algorithm finds the decision vector for a given point in the Pareto front.

\begin{figure}
	\centering
	\subfloat[F]{{ \includegraphics[width=2in]{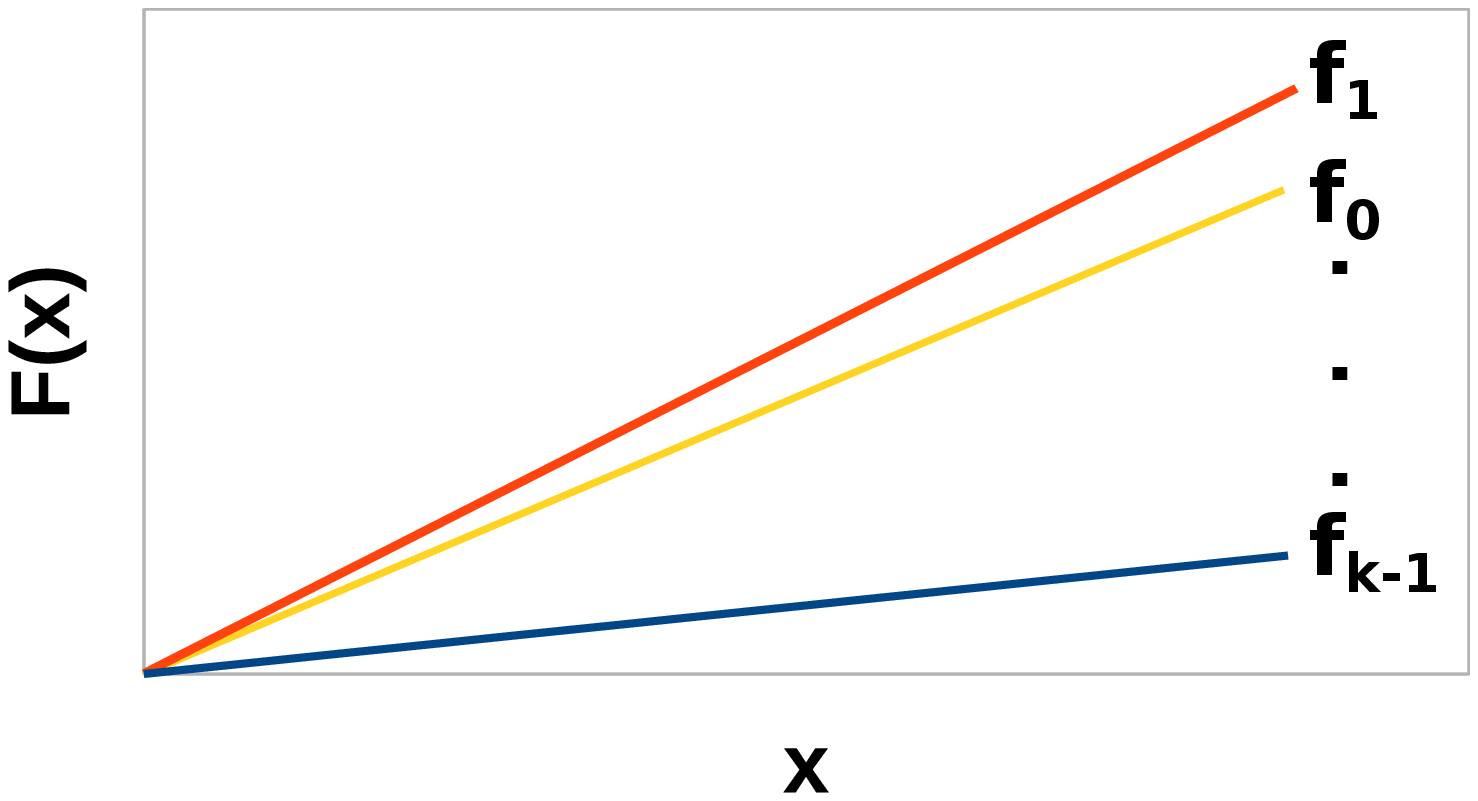} } \label{fig:linearFuncs_f}}
	\subfloat[G]{{ \includegraphics[width=2in]{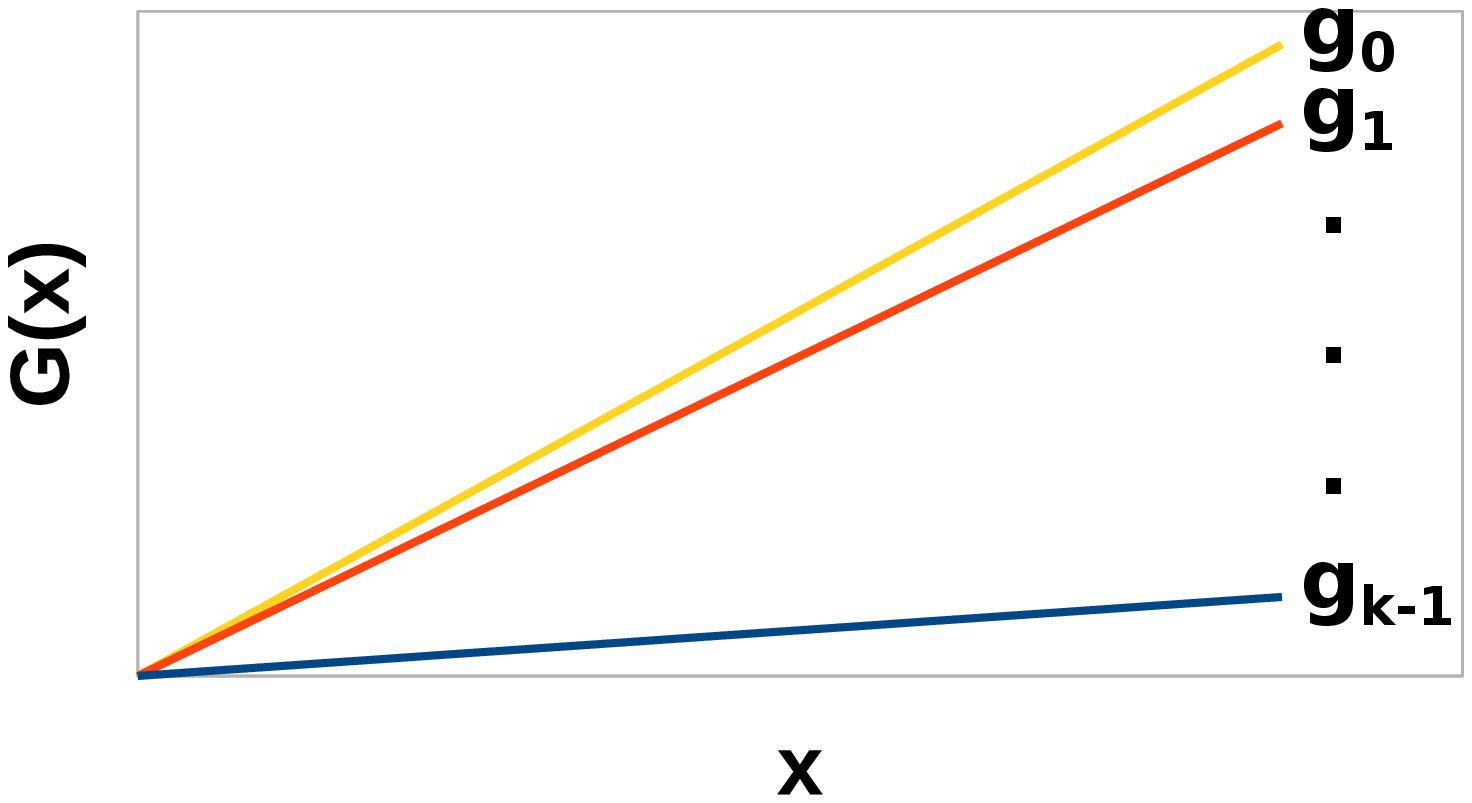} } \label{fig:linearFuncs_g}}
	\hfill
	\subfloat[Pareto front]{{ \includegraphics[width=2.5in]{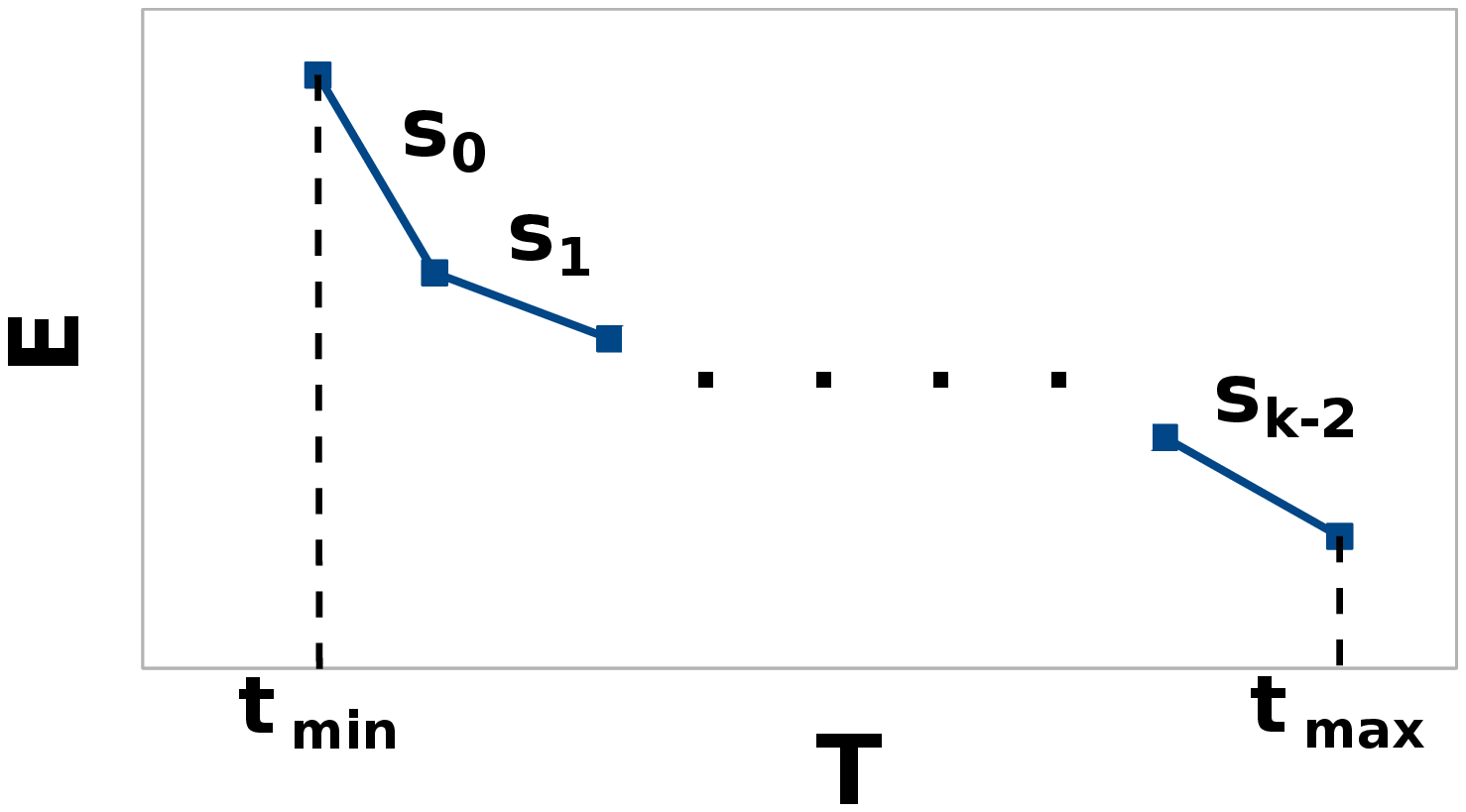} } \label{fig:linear_pareto_sample}}
	\caption{Sets F and G of $k$ linear increasing functions each. Functions in $G$ are arranged in the decreasing order of slopes. LBOPA returns a $k$-chain linear piece-wise  Pareto front shown in Fig.\ref{fig:linear_pareto_sample}.}
	\label{fig:linear_funcs_sols}
\end{figure}

The inputs to LBOPA (see Algorithm \ref{alg:lbopa} for pseudo-code) are two sets of $k$ functions each, $F$ and $G$, and an input value, $n$ $\in \mathbb{R}_{> 0}$. LBOPA constructs a Pareto front, consisting of $k-1$  segments $\{s_0, s_1, \cdots, s_{k-2}\}$. Each segment $s_i$ has two endpoints, $(t_i, e_i)$ and $(t_{i+1}, e_{i+1})$, which are connected by curve $P_f(t) =  b_i \times n - \sum_{j= i + 1}^{k - 1} (b_i - b_j)\times f_j^{-1}(t)$  ($0 \le i \le k - 2$). Figures \ref{fig:linearFuncs_f} and \ref{fig:linearFuncs_g} illustrate the functions in the sets, $F$ and $G$, when all functions in $F$ are linear, $f_i(x)=a_i \times x$. In this particular case, the Pareto front returned by LBOPA  will be piece-wise linear, $P_f(t) =  b_i \times n - t \times \sum_{j= i + 1}^{k - 1} \frac{b_i - b_j}{a_j}$ ($0 \le i \le k - 2$), as shown in Figure \ref{fig:linear_pareto_sample}.

\begin{algorithm}
	\caption{Algorithm LBOPA constructing the Pareto front of the optimal solutions,  minimizing the max of continuous and strictly increasing functions and the sum of linear increasing functions, in the objective space  $T \times E$.} 
	\label{alg:lbopa}
	\begin{algorithmic}[1]
		\Function{LBOPA}{$n,k,F,G$}
			\State $S \gets \varnothing$
			\For{$i \gets 0, k-1$}
				\State $t_i$ $\gets$ $\min_{X} \text{  } \{ \text{  } \max_{j=i}^{k-1}~f_j(x_j) \text{  } \}$
				\State $e_i$ $\gets$ $ b_i \times n - \sum_{j= i + 1}^{k - 1} (b_i - b_j)\times f_j^{-1}(t_i)$
				\State $S \gets S \cup (t_i, e_i)$
			\EndFor
			\For{$i \gets 0, k-2$}
				\State Connect $(t_i, e_i)$ and $(t_{i+1}, e_{i+1})$ by curve $b_i \times n - \sum_{j= i + 1}^{k - 1} (b_i - b_j)\times f_j^{-1}(t)$
			\EndFor
		\EndFunction
	\end{algorithmic}
\end{algorithm}

The main loop of the Algorithm \ref{alg:lbopa} computes $k$ points (Lines 3-7). In an iteration $i$, the minimum value of objective $T$, $t_i$, is obtained using the algorithm, solving the single-objective min-max optimization problem, $\min_{X} \{ \max_{j=i}^{k-1}~f_j(x_j)\}$. We do not present the details of this algorithm. Depending on the shapes of functions, $\{f_0,\dots,f_{k-1}\}$, one of the existing polynomial algorithms solving this problem can be employed \cite{Lastovetsky2004},\cite{Lastovetsky2007}.

The end point $(t_{min}, e_{max}) = (t_0,e_0)$  represents decision vectors with the minimum value of objective $T$ and the maximum value of objective $E$, while the end point  $(t_{max}, e_{min}) = (t_{k-1},e_{k-1})$ represents decision vectors with the maximum value of objective $T$ and the minimum value of objective $E$ (as illustrated for the case of all linear increasing functions in Figure \ref{fig:linear_pareto_sample}).

Given an input $t \in [t_0,t_{k-1}]$, Algorithm \ref{alg:partition} finds a decision vector $X$ = \{$x_0$, $x_1$, $\cdots$, $x_{k-1}$\} such that $\sum_{i=0}^{k-1} x_i = n$,  $\max_{i=0}^{k-1}~f_i(x_i) = t$,  and  $\sum_{i=0}^{k-1} g_i(x_i)$ is minimal. The algorithm first initialises $X$ with $\{x_0, x_1, \cdots, x_{k-1}~|~x_i= f_i^{-1}(t)\}$ (Line \ref{alg:partition:initSol}) so that $f_i(x_i)=t$ for all $i \in [0,k-1]$. For this initial $X$  the condition  $\max_{i=0}^{k-1}~f_i(x_i) = t$ is already satisfied but $\sum_{i=0}^{k-1} x_i$ may be either equal to $n$ or greater than $n$. If  $\sum_{i=0}^{k-1} x_i = n$, then this initial $X$ will be the only decision vector such that $\sum_{i=0}^{k-1} x_i = n$ and $\max_{i=0}^{k-1}~f_i(x_i) = t$ and hence the unique (Pareto-optimal) solution. Otherwise, $\sum_{i=0}^{k-1} x_i = n + n_{plus}$ where $n_{plus}>0$. In that case, this initial vector$X$ will maximize both $\sum_{i=0}^{k-1} x_i$ and $\sum_{i=0}^{k-1} g_i(x_i)$ in the set $\mathcal{X}_t$ of all vectors in the decision space satisfying the condition $\max_{i=0}^{k-1}~f_i(x_i) = t$. The algorithm then iteratively reduces elements of vector $X$ until their sum becomes equal to $n$. Obviously, each such reduction will also reduce $\sum_{i=0}^{k-1} g_i(x_i)$. To achieve the maximum reduction of $\sum_{i=0}^{k-1} g_i(x_i)$, the algorithm starts from vector element $x_i$, the reduction of which by an arbitrary amount $\Delta x$ will result in the maximum reduction of $\sum_{i=0}^{k-1} g_i(x_i)$. In our case, it will be $x_0$ as the functions in $G$ are sorted in the decreasing order of coefficients $b_i$.  Thus, at the first reduction step, the algorithm will try to reduce $x_0$ by $n_{plus}$. If $x_0 \ge n_{plus}$, it will succeed and find a Pareto-optimal decision vector $X = \{x_0 - n_{plus}, x_1, \cdots, x_{k-1}\}$. If $x_0 < n_{plus}$, it will reduce $n_{plus}$ by $x_0$, set $x_0 = 0$ and move to the second step. At the second step, it will  try to reduce $x_1$ by the reduced $n_{plus}$, and so on. This way the algorithm minimizes  $\sum_{i=0}^{k-1} g_i(x_i)$, preserving $\max_{i=0}^{k-1}~f_i(x_i) = t$ and achieving   $\sum_{i=0}^{k-1} x_i = n$.

\begin{algorithm}
\caption{Algorithm finding a Pareto-optimal decision vector $X = \{x_0, x_1, \cdots, x_{k-1}\}$ for the problem $BOPGVEC(n, k, F, G)$, where functions in $F$ are continuous and strictly increasing and functions in $G$ are linear increasing, for a given point  $(t,e)$ from the Pareto front of this problem, $(t,e) \in \Psi\downarrow_{T \times E}$. Only the first coordinate of the input point, $t$, is required for this algorithm.}
\label{alg:partition}
\begin{algorithmic}[1]
\Function{Partition}{$n,k,F,G,t$}
	\State $X = \{x_0, \cdots, x_{k-1}~|~x_i \gets f_i^{-1}(t)\}$	\label{alg:partition:initSol}
	\State $n_{plus} \gets \sum_{i=0}^{k-1} x_i - n$
	\If{$n_{plus}<0$}
		\State \Return $(0, 0, \varnothing)$
	\EndIf
	\State $i$ $\gets$ $0$
	 \While{($n_{plus} > 0) \land (i<k-1)$}									\label{alg:partition:loopS}
		\If{$x_i \ge n_{plus}$}
		    \State $x_i \gets x_i - n_{plus}$
		    \State $n_{plus} \gets 0$
		\Else
		    \State $n_{plus} \gets n_{plus} - x_i$
		    \State $x_i \gets 0$	
		    \State $i$ $\gets$ $i + 1$           
		\EndIf
	 \EndWhile													\label{alg:partition:loopE}
	\If{$n_{plus}>0$}
		\State \Return $(0, 0, \varnothing)$
	\EndIf
	\State $e \gets \sum_{i=0}^{k-1} b_i \times x_i$	\label{alg:partition:findE}
	\State \Return $(t, e, X)$
\EndFunction
\end{algorithmic}
\end{algorithm}

The correctness of these algorithms is proved in Theorem \ref{thr_alg_correct}.

\begin{theorem}	\label{thr_alg_correct}	
	Consider bi-objective optimization problem $BOPGVEC(n, k,  F, G)$ where all functions in $F$ are continuous and strictly increasing and $G = \{g_i(x)~|~g_i(x) = b_i \times x, b_i \in \mathbb R_{> 0}, i \in \{0,\cdots,k-1\}\}$. Then, the piece-wise  function $S$, returned by \Call{LBOPA}{$n,k,F,G$} (Algorithm \ref{alg:lbopa}) and consisting of $k - 1$  segments, is the Pareto front of this problem, $\Psi\downarrow_{T \times E}$, and for any $(t,e) \in \Psi\downarrow_{T \times E}$, Algorithm \ref{alg:partition} returns a Pareto-optimal  decision vector $X$ such that $T(X)=t$ and $E(X)=e$.
\end{theorem}
\textit{Proof.} First, consider Algorithm \ref{alg:partition} and arbitrary input parameters $n>0$ and $t>0$.  If after initialization of $X$ (Line \ref{alg:partition:initSol}) we will have $\sum_{i=0}^{k-1} x_i < n$, it means that $t$ is too small for the given $n$, and for any vector  $Y = \{y_0, y_1, \cdots, y_{k-1}\}$ such that $\sum_{i=0}^{k-1} y_i = n$, $\max_{i=0}^{k-1}~f_i(y_i) > t$. In this case, there is no solution to the optimization problem, and the algorithm terminates abnormally.

Otherwise, the algorithm enters the $while$ loop (Line \ref{alg:partition:loopS}). If $i<k-1$  upon exit from this loop, then the elements of vector $X$ will be calculated as 
\begin{equation}\label{eq:decision_vector}
	\begin{split}
		x_{j} = \begin{cases}	
			0 							& \quad j < i \\
			n - \sum_{m=j+1}^{k-1}  f_m^{-1}(t)		& \quad j = i \\
			f_j^{-1}(t)						& \quad j > i
		\end{cases}
	\end{split}
\end{equation}
\noindent and therefore satisfy the conditions $\sum_{j=0}^{k-1} x_j = n$ and $\max_{j=0}^{k-1}~f_j(x_j) = t$.
Moreover,  the total amount of $n$ will be distributed in $X$ between vector elements with higher indices, which have lower G cost, $g_i(x)$, because $b_{i} \ge b_{i+1}, \forall i \in \{0,\cdots,k-2\}$. Therefore, for any other vector $Y = \{y_0, y_1, \cdots, y_{k-1}\}$ satisfying these two conditions, we will have $\sum_{i=0}^{k-1} g_i(y_i) \ge \sum_{i=0}^{k-1} g_i(x_i)$. Indeed, such a vector $Y$ can be obtained from $X$ by relocating certain amounts from vector elements with higher indices to vector elements with lower indices, which will increase the G cost of the relocated amounts. Thus, when the algorithm exits from the $while$ loop with $i<k-1$, it will return a Pareto-optimal decision vector $X$. 

If  the algorithm exits from the $while$ loop with $i=k-1$, it will mean that $t$ is too big for the given $n$. We would still have $n_{plus}>0$ to take off the last vector element, $x_{k-1}$, but if we did it, we would make $\max_{j=0}^{k-1}~f_j(x_j) < t$. This way we would construct for the given $n$ a decision vector, which minimizes $\sum_{i=0}^{k-1} g_i(x_i)$  but whose $\max_{j=0}^{k-1}~f_j(x_j)$ will be less than $t$, which means that no decision vector $X$ such that $\max_{j=0}^{k-1}~f_j(x_j) = t$ can be Pareto optimal. Therefore, in this case the algorithm also terminates abnormally.

Thus, for any $t \in T$ Algorithm \ref{alg:partition} either finds a Pareto-optimal decision vector $X$ such that $T(X)=t$ and $E(X)=\sum_{i=0}^{k-1} b_i \times x_i=e$,  or returns abnormally if such a vector does not exist.  Let Algorithm \ref{alg:partition} return normally, and the loop variable $i$ be equal to $s$ upon exit from the loop. 
Then, according to formula \ref{eq:decision_vector},  $e = \sum_{i=0}^{k-1} b_i \times x_i= b_s \times (n - \sum_{i=s+1}^{k-1}  f_i^{-1}(t)) + \sum_{i=s+1}^{k-1} (b_i \times  f_i^{-1}(t)) =   b_s \times n - \sum_{i= s + 1}^{k - 1} (b_s - b_i)\times f_i^{-1}(t)$, where $s, n, b_i, b_s, a_i$ are all known constants. Therefore, the Pareto front $e = P_f(t)$ can be expressed as follows: 
\begin{equation*}	\label{eq:linearsegment}
	\begin{split}
		&e = P_f(t) = b_s \times n - \sum_{i= s + 1}^{k - 1} (b_s - b_i)\times f_i^{-1}(t)\\
		&t_{min} = \min_{X} \text{  } \{ \text{  } \max_{j=i}^{k-1}~f_j(x_j) \text{  } \}, t_{max}=  f_{k-1}(n) \\
		&t \in [t_{min}, t_{max}] \\
		&s \in \mathbb{Z}_{[0, k -2]},
	\end{split}
\end{equation*}
\noindent which is the analytical expression of the piece-wise function constructed by Algorithm \ref{alg:lbopa} (LBOPA).

\noindent \textit{End of Proof}.

\begin{theorem}	\label{thr:LinearTime}
LBOPA (Algorithm \ref{alg:lbopa}) and PARTITION (Algorithm \ref{alg:partition}) have polynomial time complexities.
\end{theorem}

\textit{Proof.} The \emph{for} loop in LBOPA (Algorithm \ref{alg:lbopa}, Lines 3-7) has $k$ iterations. At each iteration $i$, the computation of $t_i$ has a time complexity of $\bigo(k^2 \times \log_2 n)$ \cite{Lastovetsky2004}, the computation of $e_i$ has a time complexity of $\bigo(k)$, and the insertion of the point in the set $\mathcal{S}$ has complexity $\bigo(1)$. Therefore, the time complexity of the loop is $\bigo(k^2 \times \log_2 n)$. The time complexity of the loop (Lines 8-10) is $\bigo(k)$. Therefore, the time complexity of the Algorithm \ref{alg:lbopa} is $\bigo(k^2 \times \log_2 n)$.

Let us consider the PARTITION algorithm \ref{alg:partition}. The initialization of $X$ (Line 2) and computation of $n_{plus}$ has time complexity $\bigo(k)$ each. The while loop (Lines 8-17) iterates as long as $n_{plus} > 0$ and $i < k-1$, of which $i < k-1$ is the worst case scenario. The time complexity of the loop is, therefore, $\bigo(k)$. The time complexity of computation of $e$ in Line 21 is $\bigo(k)$. Therefore, the time complexity of the Algorithm \ref{alg:partition} is bounded by $\bigo(k)$.

\textit{End of Proof}.

\section{Bi-objective Optimization Problem for Max and Sum of Discrete Functions} \label{sec:gbopa}

In this section, we solve a version of BOPGVEC, called BOPGVECD, where the input $n$ is a positive integer and all the functions in $F$ and $G$ are discrete functions represented by sets of cardinality $m$. We start with the formulation of BOPGVECD. We then propose an algorithm, GBOPA, solving BOPGVECD when the functions are of an arbitrary shape. GBOPA is an exact algorithm that employs the branch-and-bound technique. To shrink the search space, it applies two bounding criteria, \emph{sum threshold} and \emph{size threshold}.

\subsection{Formulation of the Bi-objective Optimization Problem} \label{sec:bop_formulation_2}

Given a positive integer $n \in \mathbb{Z}_{> 0}$ and two sets of $k$ functions each, $F = \{f_0, f_1, \cdots, f_{k-1}\}$ and  $G = \{g_0, g_1, \cdots, g_{k-1}\}$, where $f_i, g_i \colon \mathbb{R}_{\ge 0} \to \mathbb{R}_{\ge 0}, i \in \{0,\cdots,k-1\}$ are discrete functions represented by sets of cardinality $m$, the problem is to find a vector $X = \{x_0,\cdots,x_{k-1}\} \in  \mathbb{Z}^{k}_{\ge 0}$ such that $\sum_{i=0}^{k-1} x_i = n$, minimizing the objective functions $T(X) = \max_{i=0}^{k-1}~f_i(x_i)$ and $ E(X) =  \sum_{i=0}^{k-1} g_i(x_i)$.

Thus, the problem can be formulated as follows:

\textbf{BOPGVECD($n, k, F, G$)}:
\begin{equation} \label{eq:bopfg2}
\begin{split}
T(X) &= \max_{i=0}^{k-1} f_i(x_i) \\
E(X) &= \sum_{i=0}^{k-1} g_i(x_i) \\
\underset{X}{\text{minimize}} & \quad \{T(X), E(X)\} \\
\quad \text{s.t.} & \quad x_0 + \cdots + x_{k-1} = n
\end{split}
\end{equation}

Our proposed solution, GBOPA, solves BOPGVECD by finding the discrete Pareto front containing the Pareto-optimal objective vectors and the decision vector associated with each such objective vector. GBOPA finds a set of triplets $\{(T(X),E(X),X)\}$ such that $X$ is the decision vector corresponding to the Pareto-optimal objective vector, $(T(X),E(X))$.

To aid the exposition of GBOPA, we use the following notation:
\begin{equation*}
\begin{split}
MAX(k,F,X) &= \max_{i=0}^{k-1}~f_i(x_i) \\
SUM(k,G,X) &= \sum_{i=0}^{k-1} g_i(x_i) \\
SUBMAX(i,k,F,X) &= \max_{j=i}^{k-1}~f_j(x_j)\\
SUBSUM(i,k,G,X) &= \sum_{j=i}^{k-1} g_j(x_j)
\end{split}
\end{equation*}

\subsection{Formal Description of GBOPA}	\label{sec:gbopa_GBOPAformal}

\begin{algorithm}
	\scriptsize
	\caption{Algorithm solving BOPGVECD.} \label{alg_code_GBOPA}
	\begin{algorithmic}[1]	
		\Function{GBOPA}{$n, k, F, G, \Psi$}
		\Statex \textbf{INPUT:}
		\Statex Problem size, $n \in \mathbb Z_{> 0}$
		\Statex Dimension of $X$, $k \in \mathbb Z_{> 0}$
		\Statex Constraint functions $F$, $f = \{f_0(x),...,f_{p - 1}(x)\}$,
		\Statex $f_i(x) = \{(x_{ij},f_{ij})~|~i \in [0\isep k), j \in [0 \isep m), f_{ij} \in \mathbb R_{\ge 0}\}, x_{ij} \in \mathbb Z_{> 0}$.
		\Statex Constraint functions $G$, $G = \{g_0(x),...,g_{k - 1}(x)\}$,
		\Statex $g_i(x) = \{(x_{ij},t_{ij})~|~i \in [0 \isep k), j \in [0 \isep m), g_{ij} \in \mathbb R_{\ge 0}\}, x_{ij} \in \mathbb Z_{> 0}$.
		\Statex \textbf{OUTPUT:}
		\Statex Pareto optimal solutions for the two objectives $t$ and $e$, $\Psi$,
		\Statex $\Psi = \{(t_i,e_i,X)~|~i \in [0 \isep |\Psi|)\}$,
		\Statex $X = \{x[0],x[1],\cdots,x[k-1]\}$,
		\Statex $x[i] \in \{D_i \cup \{0\}\},~i \in [0 \isep k)$.
		\Statex
		\State $F$ $\gets$ $F \cup Sort_\uparrow(F)$ $, $ $G$ $\gets$ $G \cup Sort_\uparrow(G)$	\label{gbopa_sorting1}
		
		\State $(X_{t_{min}}, t_{min})$ $\gets$ \Call {HPOPTA}{$n, k, F$}		\label{gbopa_optE_time}
		\State $\varepsilon$ $\gets$ $SUM(k,G,X_{t_{min}})$	\label{gbopa_e-threshold}
		\State $\sigma$ $\gets$ \Call {SizeThresholdCalc}{$k,G,\varepsilon$}		\label{gbopa_s_th}
		\State $PMem[i][j]$ $\gets$ $\emptyset$, $\forall i \in \{1,\cdots,k-2\},$ $j \in \{0,\cdots,n\}$ \label{gbopa_mem}
		\State \Call {GBOPA\_Kernel}{$0, n, k, F, G, \varepsilon, \sigma, X_{cur}, PMem, \Psi$} \label{gbopa_kernelCall}
		\State \Return $\Psi$
		\EndFunction		
	\end{algorithmic}
\end{algorithm}

\begin{algorithm}
	\scriptsize
	\caption{The core recursive kernel invoked by GBOPA.} \label{hep_kernel_code}
	\begin{algorithmic}[1]	
		\Function{GBOPA\_Kernel}{$lvl, n, k, F, G, \varepsilon, \sigma, X_{cur}, PMem, \Psi$}
		\Statex
			\If{\Call{Cut}{$n, \sigma_{lvl}$}}					\label{hep_kernel_size_thr1}
				\State \Return $FALSE$
			\EndIf												\label{hep_kernel_size_thr2}
			\If{$lvl = k - 1$ $\wedge$ $g_{lvl}(n) \le \varepsilon$}		\label{hep_kernel_leaf1}
				\State $x_{cur}[lvl]$ $\gets$ $n$
				\State \Return $TRUE$
			\Else
				\State \Return $FALSE$
			\EndIf									\label{hep_kernel_leaf2}
			\If{$n \neq 0$ $\wedge$ $lvl \ge 1$ $\wedge$ $lvl \leq k - 2$}	\label{hep_kernel_retrieveMem1}
				\State $status$ $\gets$ \Call{ReadParetoMem}{$n, lvl, \varepsilon, PMem$} \label{hep_kernel_rSolMem}
				\If{$status=SOLUTION$}
  				    \State \Return $TRUE$
				\EndIf				
				\If{$status=NOT\_SOLUTION$}
					\State \Return $FALSE$
				\EndIf
			\EndIf									\label{hep_kernel_retrieveMem2}
			\State $isSol$ $\gets$ $FALSE$			\label{hep_kernel_setzero}
			\State $partsVec \gets \varnothing$			
			\State $indx$ $\gets$ $-1$
			\State $x_{lvl~indx}$ $\gets$ $0$
			\While{$g_{lvl}(x_{lvl~indx})  \le \varepsilon$} 			\label{hep_kernel_mainLoop1}
				\If{$x_{lvl~indx} \le n$}			
					\State $x_{cur}[lvl]$ $\gets$ $x_{lvl~indx}$	\label{hep_kernel_store_dc}									
					\State $outRes$ $\gets$ \Call{GBOPA\_Kernel}{$lvl+1,n-x_{lvl~indx},k,F, G, \varepsilon,\sigma,x_{cur}, PMem, \Psi$}	\label{hep_kernel_recall}
					\If{$outRes = TRUE$}	\label{hep_kernel_onreturn1}
						\State $isSol$ $\gets$ $TRUE$
						\State $partsVec$ $\gets$ $partsVec$ $\cup$ $x_{lvl~indx}$
					\EndIf					\label{hep_kernel_onreturn2}
				\EndIf				
				\If{$n = 0$ $\vee$ $indx+1=m$}				\label{hep_kernel_break1}
					\State \textbf{break}	
				\EndIf						\label{hep_kernel_break2}
				\State $indx \gets indx+1$	\label{hep_kernel_nextpoint}
			\EndWhile								\label{hep_kernel_mainLoop2}
			\If{$lvl \ge 1 \wedge  lvl \leq k - 2$}	\label{hep_kernel_merge1}
				\State \Call {MergePartialParetoes}{$n, k, lvl, F, G, partsVec, PMem,\Psi$}	
			\EndIf									\label{hep_kernel_merge2}
			\State \Call{MakeParetoFinal}{$PMem[lvl][n]$}	\label{hep_kernel_fin}
			\State \Return $isSol$	\label{hep_kernel_return}			
		\EndFunction		
	\end{algorithmic}
\end{algorithm}

GBOPA is illustrated by the Algorithm \ref{alg_code_GBOPA}. Its inputs are: the input problem size, $n$; the sets of discrete functions, $F=\{f_0(x),\cdots,f_{k-1}(x)\}$ and $G = \{g_0(x),\cdots,g_{k-1}(x)\}$. Each function, $f_i(x)$ in $F$ is represented by an array of $m$ pairs, $(x_{ij},f_{ij})$, $j \in \{0,1,\cdots,m-1\}$, so that $x_{ij}$ is the $j$-th data point in the function and $f_{ij}$ represents $f_i(x_{ij})$. Similarly, each function $g_i(x) \in G$ is represented by $m$ pairs $(x_{ij},g_{ij})$, $j \in \{0,1,\cdots,m-1\}$, so that $x_{ij}$ is the $j$-th data point in the function and $g_{ij}$ represents $g_i(x_{ij})$. GBOPA returns $\Psi$, the set of Pareto optimal solutions. It consists of triples, $(MAX(k,F,X), SUM(k,G,X), X)$ where $X=\{x_0,\cdots,x_{k-1}\}$ is a Pareto optimal solution and $MAX(k,F,X)$ and $SUM(k,G,X)$ are the corresponding objective values. The solutions are output in increasing order of $SUM(k,G,X)$ in $\Psi$.

GBOPA starts by sorting the array $F$ ($G$) in non-decreasing values for $f_i(x)$ ($g_i(x)$), $i \in \{0,\cdots,k-1\}$ (Line \ref{gbopa_sorting1}). Both original and sorted functions are kept. Then, $\Call{HPOPTA}{}$ is invoked to solve the min-max single objective optimization problem that finds the integer partition of $n$ given the input set of functions, $F$ (Line \ref{gbopa_optE_time}). This function returns the optimal solution, $X_{t_{min}}$, and the corresponding objective value, $SUM(k,G,X_{t_{min}})$. The sum threshold $\varepsilon$ is initialized to $X_{t_{min}}$ (Line \ref{gbopa_e-threshold}). 

The size threshold array $\sigma$ is initialised by using the function \Call{SizeThresholdCalc}{} (Line \ref{gbopa_s_th}). A 2D array $PMem$, with dimensions of $(k-2) \times (n+1)$, is defined to save Pareto optimal solutions for levels $\{L_1,\cdots,L_{k-2}\}$, which are found during the tree exploration (Line \ref{gbopa_mem}). Then, \Call{GBOPA\_Kernel}{} is invoked to explore the solution tree and outputs the set of Pareto optimal solutions, $\Psi$. The structure of $PMem$ is described in the Appendix (Section \ref{app:pmem}).

An informal description of GBOPA is presented in the Appendix (Section \ref{app:gbopa}). 

\subsubsection{Recursive Algorithm GBOPA\_Kernel} \label{sec_code_heopta}

Algorithm \ref{hep_kernel_code} illustrates the core recursive function of GBOPA, \textit{GBOPA\_Kernel}. It recursively explores the solution tree and builds Pareto optimal solutions by merging the partial solutions. Pareto optimal solutions for a given node at level $L_i$, $i \in \{0,1,\cdots,k-2\}$, are built by merging all solutions stored for its children, placed at level $L_{i + 1}$. To reduce the search space and thus to achieve a polynomial computational complexity, \textit{GBOPA\_Kernel} uses three key operations \emph{Cut}, \emph{SavePareto} and \emph{ReadParetoMem}, described in the section \ref{sec:gbopa}. 

The input variable $lvl$ indicates the tree level that is processing in the current recursion of \textit{GBOPA\_Kernel}. Prior to expanding a node at the level $L_{lvl}$, \textit{GBOPA\_Kernel} determines whether its workload exceeds its correspondence bounding criterion $\sigma_{lvl}$. If it is the case then the node is not explored (Lines \ref{hep_kernel_size_thr1}-\ref{hep_kernel_size_thr2}). Lines \ref{hep_kernel_leaf1}-\ref{hep_kernel_leaf2} process solutions found at the last level $L_{k-1}$ in the tree. If there exists a solution, the function returns $TRUE$, otherwise $FALSE$.

Before exploring a node at a given level $lvl$, $lvl \in \{1,2,\cdots,k-2\}$, the function \Call{ReadParetoMem}{} is called to retrieve from  the partial solutions previously saved for the current workload $n$ on the level $L_{lvl}$ (Lines \ref{hep_kernel_retrieveMem1}-\ref{hep_kernel_retrieveMem2}). The variable \emph{status} determines the type of solutions retrieved from $PMem$. If no solution is already stored for the node, \emph{GBOPA\_Kernel} returns $FALSE$ and backtracks.  If $e' \ge \varepsilon$ for all the retrieved signified by the status, \emph{NOT\_SOLUTION}, \emph{GBOPA\_Kernel} returns $FALSE$ and backtracks.

If at least one of the partial solutions, in the set retried, has a $e'$ value less than $\varepsilon$ (status, \emph{SOLUTION}), the function returns $TRUE$. If none of the above cases happen, the routine starts expanding the node by initializing pointer $indx$ to $-1$ and $x_{lvl~indx}$ to $0$ (Lines \ref{hep_kernel_setzero}-\ref{hep_kernel_mainLoop2}). The variable $indx$, ranging from $-1$ to $m-1$, determines the data points in the discrete functions where $x_{lvl~indx}$ represents the value of $indx$-th data point in the $g_{lvl}(x)$ ($f_{lvl}(x)$).

The $while$ loop (Lines \ref{hep_kernel_mainLoop1}-\ref{hep_kernel_mainLoop2}) examines all data points $x_{lvl~indx} \in D_{lvl}$ where $g_{lvl}(x_{lvl~indx})$ is less than or equal to $\varepsilon$. The array $X_{cur}=\{x_{cur}[0],\cdots,x_{cur}[p-1]\}$ stores data points currently assigned to $x_i, i \in \{0,1,\cdots,k-1\}$. In each iteration, the data point $x_{lvl~indx}$ is extracted from $D_{lvl}$ and stored in array $x_{cur}[lvl]$ (Line \ref{hep_kernel_store_dc}). \emph{GBOPA\_Kernel} is recursively invoked to find solutions for the remaining problem size $n-x_{lvl~indx}$ at the next level $L_{lvl+1}$ (Line \ref{hep_kernel_recall}). If there exists any solution for the workload, $x_{lvl~indx}$ is added to $partsVec$, a list holding all data points which result in Pareto optimal solutions on level $L_{lvl}$ (Lines \ref{hep_kernel_onreturn1}-\ref{hep_kernel_onreturn2}). The variable $indx$ is incremented to examine all data point in $D_{lvl}$ one-by-one, and the \emph{while} loop terminates whenever all the data points are examined (Lines \ref{hep_kernel_break1}-\ref{hep_kernel_break2}).

After exploring all children of the current node, the function \Call{MergePartialParetoes}{} is invoked to build and store the partial solutions of the node by merging the partial solutions of its children.

In the end, the corresponding array cell storing the Pareto optimal solution for a node with a problem size $n$ at $L_{lvl}$ ($PMem[lvl][n]$) is labelled \emph{Finalized} (Line \ref{hep_kernel_fin}). Finalizing a memory cell implies that this cell contains the final partial solutions. \textit{GBOPA\_Kernel} returns $TRUE$ provided that exploring the node, processed in the current recursion, leads to a solution (Line \ref{hep_kernel_return}).

The pseudocodes of the following subroutines are described in the Appendix (Section \ref{app:subroutines}). 
\begin{itemize}
	\item The function \emph{SizeThresholdCalc}, which calculates the size threshold array, $\sigma$.
	\item The function \emph{Cut}, which checks if the input workload $n$ is greater than the input size threshold $\sigma$.
	\item The function \emph{ReadParetoMem} that retrieve the saved solutions for a given workload $n$ on Level $L_{lvl}$.
	\item The function \emph{MakeParetoFinal}, which finalizes the input memory cell, $pmem$.
	\item The function \emph{MergePartialParetoes}, which builds the Pareto front of a node using the Pareto fronts of its children.
\end{itemize}

\subsection{Correctness Proof of GBOPA}

\begin{theorem}
Consider BOPGVECD where each function in the sets, $F$ and $G$, is represented by a set of cardinality $m$. The algorithm GBOPA solves BOPGVECD and returns the set of Pareto optimal solutions.
\end{theorem}

\textit{Proof.} A naive algorithm to solve the problem builds the full tree of solutions and then determines the Pareto front. However, it has exponential complexity. GBOPA achieves polynomial complexity by using memorization and applying two bounding criteria employed in the \emph{Cut} operation, which explores only a small fraction of the full solution tree to determine the Pareto front. Therefore, the correctness of GBOPA will be proved if we show that there exists no subtree that contains a Pareto optimal solution and ignored by \emph{Cut}.

Consider a solution, $X=\{x_0,\cdots,x_{k-1}\}$, which is eliminated from the search space by using the \emph{Cut} operation. Using the definition of size threshold and sum threshold, $SUM(k, G, X)$ is greater than $\varepsilon$. $MAX(k, F, X)$ will also be greater than or equal to the optimal solution $X_{t_{min}}$ ($t_{min} = MAX(k,F,X_{t_{min}}), t_{min} \le MAX(k,F,X)$). As explained in Section \ref{sec:gbopa}, $\varepsilon$ is equal to $SUM(k,G,X_{t_{min}})$. Hence, there is a distribution, $X^* = \{x_0^*, x_1^*\cdots,x_{k-1}^*\}$ where its $MAX(k,F,X^*) = t_{min}$ and $SUM(k,G,X^*) = \varepsilon$. Thus, we have $MAX(k,F,X^*) < MAX(k,F,X)$ and $SUM(k,G,X^*) \le SUM(k,G,X)$, and according to the definition of Pareto optimality, the solution $X$, which is removed by \emph{Cut}, is dominated by the solution $X^*$ and cannot be a member of the Pareto optimal set.
\textit{End of Proof}.

\subsection{Complexity Proof of GBOPA}

\begin{lemma}	\label{lemma:paretoSols}
Consider BOPGVECD where each function in the sets, $F$ and $G$, is represented by a set of cardinality $m$. The maximum number of Pareto optimal solutions is equal to $k \times m$.
\end{lemma}

\textit{Proof.} We know that the objective $MAX(k,F,X)$ is a max function of $f_i(x)$, $i \in \{0, 1, \cdots, k - 1\}$, where the cardinality of each function is $m$. Therefore, there are no more than $k \times m$ unique values for the objective function. By the definition of Pareto optimality, the values for the two objectives in the Pareto front are unique. Hence, one can conclude that the maximum number of Pareto optimal solutions cannot exceed $k \times m$.
\textit{End of Proof}.

\begin{lemma}	\label{lemma:mergeComplexity}
Consider BOPGVECD where each function in the sets, $F$ and $G$, is represented by a set of cardinality $m$. The computational complexity of the function \Call{MergePartialParetoes}{} is equal to $\bigo(m^2 \times k \times \log_2(m \times k))$.   
\end{lemma}

\textit{Proof.} Consider a node $N$ at the level $L_{lvl}$ of a solution tree. As explained in Section \ref{sec:gbopa}, the node has $m + 1$ children in the general case. Using Lemma \ref{lemma:paretoSols}, each child of the node at level $L_{lvl+1}$ has up to $m \times (k-lvl-1)$ Pareto optimal solutions. Therefore, there are in total $(m + 1) \times (m \times (k-lvl-1))$ solutions, which should be examined one by one to build the Pareto optimal solutions of node $N$. 

Let $\Psi_N$ be the set storing partial Pareto optimal solutions of the node. As proved in Lemma \ref{lemma:paretoSols}, the cardinality of $\Psi_N$ does not exceed $m \times (k-lvl)$ Pareto optimal solutions. Since $\Psi_N$ is a sorted set, the cost of inserting a solution into it or removing a solution from it is logarithmic in its size, which is $\log_2(m \times (k-lvl))$. Therefore, building $\Psi_N$ has a complexity of $\bigo((m + 1) \times (m \times (k-lvl-1)) \times \log_2(m \times (k-lvl)))$ $\approxeq$ $\bigo(m^2 \times k \times \log_2(m \times k))$.

\textit{End of Proof}

\begin{lemma}	\label{lemma_1}
Consider BOPGVECD where each function in the sets, $F$ and $G$, is represented by a set of cardinality $m$. The computational complexity of \emph{GBOPA\_Kernel} is $\bigo(m^3 \times k^3 \times \log_2(m \times k))$.
\end{lemma}

\textit{Proof.} Since \emph{GBOPA\_Kernel} is a recursive algorithm, its computational complexity is determined using the number of its recursions. We will formulate the number of recursions using a sample tree depicted in the Appendix (Figure \ref{fig:search_tree_time_com}).

Consider a BOPGVECD with an input problem size $n$, five discrete functions, $F = \{f_0(x), f_1(x), f_2(x), f_3(x), f_4(x)\}$, and five discrete functions, $G$ = \{$g_0(x)$, $g_1(x)$, $g_2(x)$, $g_3(x)$, $g_4(x)$\}, each with a cardinality of 2 ($m = 2$). Suppose the domain of each function in $F$ and $G$ is the set $D = \{\Delta x, 2\Delta x\}$. It should be noted that \emph{GBOPA\_Kernel} is able to deal with functions with any step size between the points in its domain. Using a constant step size $\Delta x$ does not make the proof less general. Without loss of generality and for the sake of simplicity, we assume that for all functions in $F$ and $G$, if $x_1 < x_2$ then $\forall i \in \{0, 1, 2, 3, 4\}$, $f_i(x_1) < f_i(x_2)$ and $g_i(x_1) < g_i(x_2)$.

The figure in the Appendix (Figure \ref{fig:search_tree_time_com}) shows the solution tree exploring all possible solutions. Suppose $n$ is greater than $8\Delta x$, which is equal to the maximum possible size subtracted from $n$ in this example. In the figure, red nodes are ones that have been already expanded in the same level, and their solutions are retrieved from $PMem$. For the sake of simplicity, we only consider the two operations \emph{SavePareto} and \emph{ReadParetoMem}, and \emph{Cut} is ignored.

According to the sample tree, the number of recursions (the number of nodes whose solutions are not retrieved from the memory) in each level explored can be obtained using the Eq. \ref{eq:gbopa_recCalls}.
\begin{equation}	\label{eq:gbopa_recCalls}
    C\#(L) = \begin{cases}
    			L \times m+1	  		& 0 \leq L < k-1\\
    			C\#(k-2)\times (m + 1)	& L = k - 1\\
    		\end{cases}
\end{equation}
where $L$ represents the level number.
The expanded form of Eq. \ref{eq:gbopa_recCalls} is shown in Eq. \ref{eq:gbopa_recCalls_simple}.
\begin{equation}	\label{eq:gbopa_recCalls_simple}
    C\#(L) = \begin{cases}
               L \times m+1	  										& 0 \leq L < k-1\\
               m^2 \times k - 2 \times m^2 + m \times k - m + 1		& L = k - 1\\
           \end{cases}
\end{equation}
That is, the total number of recursive calls is equal to $\sum_{L=0}^{k-1}(C\#(L))$ which is equal to $\bigo(m \times k ^ 2 + m ^ 2 \times k)$.

In addition, Eq. \ref{eq:gbopa_mem_total} formulates the number of nodes whose results are retrieved from $PMem$ in each level.
\begin{equation}	\label{eq:gbopa_mem_total}
	\text{Memory\#(L)} = (C\#(L-1) - 1) \times m ~=~ (m^2) \times (L - 1), \quad 1 \leq L \leq k - 2\\
\end{equation}
Since $PMem$ saves the solutions which are found on levels $1$ to $k-2$, the total number of nodes whose solutions are saved (red nodes in the figure) is equal to $\sum_{L=1}^{k-2} Memory\#(L) = \bigo(m^2 \times k^2)$. Since the complexity of \Call{ReadParetoMem}{} is $\bigo(1)$, the computational cost for retrieving all solutions from $PMem$ is equal to $\bigo(m^2 \times k^2)$.

The function \Call{MergePartialParetoes}{} is invoked after exploring all children of any node (black nodes in Figure \ref{fig:search_tree} in the Appendix) in levels $\{L_0,\cdots,L_{k-2}\}$. Regarding Lemma \ref{lemma:mergeComplexity} and Eq. \ref{eq:gbopa_recCalls}, the total cost of all \Call{MergePartialParetoes}{} calls is equal to $\sum_{L=0}^{k-2} (L \times m + 1) \times (m^2 \times k \times \log_2(m \times k)) = \bigo(m^3 \times k^3 \times \log_2(m \times k))$.

The computational complexity of \emph{GBOPA\_Kernel} can therefore be summarized as follows:
\begin{equation*}
	\begin{split}
		\text{Complexity(\emph{GBOPA\_Kernel})} = & \bigo(\text{recursive calls of \emph{GBOPA\_Kernel}}) + \\
		& \bigo(\text{$PMem$ solutions}) + \\
		& \bigo(\text{\Call{MergePartialParetoes}{} calls}).
	\end{split}
\end{equation*}
which equals:
\begin{equation*}
	\begin{split}
		\text{Complexity(\emph{GBOPA\_Kernel})} = & \bigo(m \times k^2 + m^2 \times k)+ \\
		& \bigo(m^2 \times k^2) +\\
		& \bigo(m^3 \times k^3 \times \log_2(m \times k)) \\
		& = \bigo(m^3 \times k^3 \times \log_2(m \times k)).
	\end{split}
\end{equation*}

\begin{theorem}
Consider BOPGVECD where each function in the sets, $F$ and $G$, is represented by a set of cardinality $m$. The computational complexity of GBOPA is $\bigo(m^3 \times k^3 \times \log_2(m \times k))$.
\end{theorem}

\textit{Proof.} GBOPA consists of following main steps:
\begin{itemize}
	\item \textbf{Sorting:} There exist $k$ discrete functions in $F$ and $k$ discrete functions in $G$, each with a cardinality of $m$. The complexity to sort all the functions is $\bigo(k \times m \times \log_2 m)$.
	\item \textbf{Initializing sum threshold $\varepsilon$:} Obtaining the sum threshold involves two steps: (i). Invoking \emph{HPOPTA}, which has the complexity of $O (m^3 \times k^3)$ \cite{khaleghzadeh2018novel}, and (ii). Calculating the sum threshold $\varepsilon$, which has the complexity of $\bigo(k)$. Therefore, the complexity of this step is equal to $O (m^3 \times k^3)$.
	\item \textbf{Finding size thresholds:} To find the size threshold of a given level $L_i$, $i \in \mathbb{Z}_{[0,k-1]}$, in the worst case, all data points existing in the domain of $g(x)$ should be examined, which has a complexity of $\bigo(m)$. Therefore, finding $k$ size thresholds costs $\bigo(k \times m)$.
	\item \textbf{Memory initialization:} $PMem$ is a $(k-2) \times (n+1)$ matrix. In this step, all cells in $PMem$ are initialized, which costs $\bigo(k \times n)$.
	\item \textbf{Kernel invocation:} According to Lemma \ref{lemma_1}, the complexity of GBOPA\_Kernel is $\bigo(m^3 \times k^3 \times \log_2(m \times k))$.
\end{itemize}

Thus, the computational complexity of GBOPA is equal to the summation of all these steps, which is equal to $\bigo(m^3 \times k^3 \times \log_2(m \times k))$.
\textit{End of Proof}.

\begin{theorem}	\label{prop_mem_use}
Consider BOPGVECD where each function in the sets, $F$ and $G$, is represented by a set of cardinality $m$. The total memory consumption of GBOPA is $\bigo(n \times m \times k^2)$.
\end{theorem}

\textit{Proof.} GBOPA memories the following information:
\begin{itemize}
	\item \textbf{F functions:} There are $k$ discrete functions in $F$, each with a cardinality of $m$. We store two sorted sets, one sorted by the domain values ($x$) and the other sorted by the range values ($f(x)$). The storage complexity is $\bigo(k \times m)$.
	\item \textbf{G functions:} There are $k$ discrete functions in $G$, each with a cardinality of $m$. We store two sorted sets, one sorted by the domain values ($x$) and the other sorted by the range values ($g(x)$). The storage complexity is $\bigo(k \times m)$
	\item \textbf{$\Psi$:} The maximum number of Pareto optimal solutions is $m \times k$ (Lemma \ref{lemma:paretoSols}). Since the dimension of each solution is equal to $k$, the maximum size of the set is $\bigo(m \times k^2)$.
	\item \textbf{PMem:} This is a matrix consisting of $(k-2) \times (n+1)$ cells. Each cell stores up to $m \times k$ Pareto optimal solutions (Lemma \ref{lemma:paretoSols}). Therefore, the memory usage of $PMem$ is equal to $\bigo(n \times m \times k^2)$.
	\item \textbf{Memory consumption of HPOPTA:} The memory usage by HPOPTA algorithm is $O (k \times (m + n))$ \cite{khaleghzadeh2018novel}.
	\item \textbf{$X_{cur}$}: This is an array of size $k$ to store the data points in the current solution.
	\item \textbf{partsVec}: This is a vector of size $\bigo(m)$ storing the data points examined at each level that result in a solution. There exist $k-1$ such vectors, one per level. The total memory consumed is, therefore, equal to $\bigo(m \times k)$.
\end{itemize}

Thus, the total memory usage of GBOPA is equal to $\bigo(n \times m \times k^2)$. 
\textit{End of Proof}.

\section{Conclusion} \label{sec:conclusion}

Bi-objective optimization problems where one objective function is \textit{max} and the other is \textit{sum} are common in the category of multiple objective minimum spanning tree problems that have important applications in the field of network design and optimization (transportation and communication networks, for example), flowshop group scheduling, after disaster blood supply chain management, collaborative production planning, and performance and energy optimization of high computing systems and applications. State-of-the-art solution methods for such problems consider objective functions and constraints that are linear functions of the decision variables. 

In this work, we considered a bi-objective optimization problem that aims to determine an optimal $k$-partition of a real number $n$  minimizing the maximum of $k$-dimensional vector of functions of objective type one and the sum of $k$-dimensional vector of functions of objective type two. This problem arises in optimization of applications for performance and energy consumption on modern heterogeneous computing platforms \cite{HamidTPDS2020}.

We proposed two algorithms. The first algorithm solved the  problem for the case where all the functions of objective type one are continuous and strictly increasing, and all the functions of objective type two are linear increasing. The second algorithm solved the problem where the input $n$ is a positive integer and all the functions are discrete functions of an arbitrary shape. Both algorithms exhibit polynomial complexity.

\section{Acknowledgments}
This publication has emanated from research conducted with the financial support of Science Foundation Ireland (SFI) under Grant Number 14/IA/2474. 

\section{Appendix}

\subsection{Introduction}

The supplementary material includes:

\begin{itemize}
    \item Solution methods for multi-objective optimization. Formulations of BOPGVEC in some well-known methods.
	\item The informal description of the algorithm GBOPA solving BOPGVECD.
	\item Formal descriptions of the subroutines in the algorithm GBOPA solving BOPGVECD.	
	\item Description of the experimental platform and scientific applications.
\end{itemize}

\subsection{Multi-Objective Optimization: Solution Methods}

There are several classifications for methods solving multi-objective optimization problems \cite{Miettinen1999},\cite{Talbi2009}. Since the set of Pareto-optimal solutions are partially ordered, one classification is based on the involvement of decision maker in the solution method to select specific solutions. This involvement can be classed as follows:
\begin{itemize}
	\item \textbf{No preference}: No preference information is provided. In this case, the solution method can return one or more solutions that the decision maker can accept or reject.
	\item \textbf{A posteriori}: The solution method determines the Pareto-optimal set of solutions and presents to the decision maker, who then selects the most preferred. These methods usually have high computational complexity. 
	\item \textbf{A priori}: The decision maker provides his or her preferences before the solution process. The most common approach is called the aggregation method where the objective functions are combined into a single objective function using the utility function, which is provided in a mathematical form by the decision maker. The single objective function is then optimized using traditional optimization methods.
	\item \textbf{Interactive}: There is repeated interaction between the decision maker and the solution method. The basic steps are \cite{Miettinen1999}:
	a). Find an initial feasible solution
	b). Interact with the decision maker
	c). Obtain a new solution. If the new solution or one of the previous solutions are acceptable to the decision maker, stop. Otherwise, go to step (b).
\end{itemize}

\subsection{Multi-Objective Optimization: A posteriori Solution Methods} \label{app:solmethods}

In this section, we show the formulations of BOPGVEC in some of the well-known \textit{A posteriori} methods.

\subsubsection{Weighting Method}

In this method, the multiple objective functions are transformed into a single objective function by associating a weighting coefficient to each objective function \cite{Miettinen1999}. BOPGVEC can be formulated as follows:
\begin{alignat*}{3}
minimize & \quad w_1 \times (\max_{i=1}^{q}~f_i(x_i)) + w_2 \times (\sum_{i=1}^{q} g_i(x_i)) & & \\
\text{Subject to} & \quad x_1 + x_2 + ... + x_q = n & & \\
& x_i \geq 0 \qquad i = 1,...,q & \\
& x_i \leq n \qquad i = 1,...,q & \\
& w_1 + w_2 = 1 & \\
& 1 \leq q \leq k \\
\text{where} \quad & k, q, n, x_i \in \mathbb Z_{> 0}, \\
\quad & f_i(x_i), g_i(x_i) \in \mathbb R_{> 0}, w_1, w_2 \in \mathbb R
\end{alignat*}

The solution of the weighting method is weakly Pareto-optimal under no additional assumptions. It is Pareto-optimal if the weighting coefficients are positive, $w_1 > 0, w_2 > 0$. The solution of the weighting method is properly Pareto-optimal if all the weighting coefficients are positive \cite{Miettinen1999}. 

Different Pareto-optimal solutions can be obtained by the weighting method by altering the positive weighting coefficients \cite{Miettinen1999}. The weakness of the weighting method is that all of the Pareto-optimal solutions cannot be found if the problem is nonconvex.

\subsubsection{$\epsilon$-Constraint Method}

In this method, one of the objective functions is chosen to be optimized and all the other objective functions are converted into constraints by setting an upper bound for each one of them \cite{Miettinen1999}. BOPGVEC can be formulated as follows:
\begin{alignat*}{3}
minimize & \quad \max_{i=1}^{q}~f_i(x_i) \\
\text{Subject to} & \quad \sum_{i=1}^{q} g_i(x_i) \leq \epsilon \\
& x_1 + x_2 + ... + x_q = n & & \\
& x_i \geq 0 \qquad i = 1,...,q & \\
& x_i \leq n \qquad i = 1,...,q & \\
& 1 \leq q \leq k \\
\text{where} \quad & k, q, n, x_i \in \mathbb Z_{> 0}, \\
\quad & f_i(x_i), g_i(x_i), \epsilon \in \mathbb R_{> 0}
\end{alignat*}

The solution of the $\epsilon$-constraint method is weakly Pareto-optimal under no additional assumptions. Theoretically, all Pareto-optimal solutions can be found by the $\epsilon$-constraint method by altering the $\epsilon$ and the function to be minimized regardless of the convexity of the problem \cite{Miettinen1999}.

Alternatively one can choose to minimize the objective of energy:
\begin{alignat*}{3}
minimize & \quad \quad \sum_{i=1}^{q} g_i(x_i) \\
\text{Subject to} & \quad \max_{i=1}^{q}~f_i(x_i) \leq \epsilon \\
& x_1 + x_2 + ... + x_q = n & & \\
& x_i \geq 0 \qquad i = 1,...,q & \\
& x_i \leq n \qquad i = 1,...,q & \\
& 1 \leq q \leq k \\
\text{where} \quad & k, q, n, x_i \in \mathbb Z_{> 0}, \\
\quad & f_i(x_i), g_i(x_i), \epsilon \in \mathbb R_{> 0}
\end{alignat*}

\subsubsection{Hybrid Method}

In this method, the basic methods, weighting method and $\epsilon$-constraint method, are combined \cite{Miettinen1999}. BOPGVEC can be formulated as follows:
\begin{alignat*}{3}
minimize & \quad w_1 \times (\max_{i=1}^{q}~f_i(x_i)) + w_2 \times (\sum_{i=1}^{q} g_i(x_i)) & & \\
\text{Subject to} & \quad \max_{i=1}^{q}~f_i(x_i) \leq \epsilon_1 \\
& c \times \max_{i=1}^{q}~f_i(x_i) + \sum_{i=1}^{q} g_i(x_i) \leq \epsilon_2 \\
& x_1 + x_2 + ... + x_q = n & & \\
& x_i \geq 0 \qquad i = 1,...,q & \\
& x_i \leq n \qquad i = 1,...,q & \\
& 1 \leq q \leq k \\
& w_1 + w_2 = 1 & \\
\text{where} \quad & k, q, n, x_i \in \mathbb Z_{> 0}, \\
\quad & f_i(x_i), g_i(x_i), \epsilon_1, \epsilon_2 \in \mathbb R_{> 0}, w_1, w_2 \in \mathbb R
\end{alignat*}

The hybrid method combines all the positive features of the basic methods and can find, theoretically speaking, all Pareto-optimal solutions regardless of the convexity of the problem.

\subsubsection{Method of Weighted Metrics}

The \textit{weighted $L_p$}-formulation of BOPGVEC follows \cite{Miettinen1999}:
\begin{alignat*}{3}
minimize & \quad (w_1 \times |\max_{i=1}^{q}~f_i(x_i) - z_1^*|^d & & \\
& + w_2 \times |\sum_{i=1}^{q} g_i(x_i) - z_2^*|^d)^d & & \\
\text{Subject to} & \quad x_1 + x_2 + ... + x_q = n & & \\
& x_i \geq 0 \qquad i = 1,...,q & \\
& x_i \leq n \qquad i = 1,...,q & \\
& w_1 + w_2 = 1 & \\
& 1 \leq q \leq p \\
\text{where} \quad & k, q, n, x_i \in \mathbb Z_{> 0}, \\
\quad & f_i(x_i), g_i(x_i), \epsilon_1, \epsilon_2 \in \mathbb R_{> 0}, w_1, w_2 \in \mathbb R_{\geq 0}
\end{alignat*}

The weighted Tchebycheff problem formulation has the form \cite{Miettinen1999}:
\begin{alignat*}{3}
minimize & \quad (w_1 \times |\max_{i=1}^{q}~f_i(x_i) - z_1^*| & & \\
& + w_2 \times |\sum_{i=1}^{q} g_i(x_i) - z_2^*|) & & \\
\text{Subject to} & \quad x_1 + x_2 + ... + x_q = n & & \\
& x_i \geq 0 \qquad i = 1,...,q & \\
& x_i \leq n \qquad i = 1,...,q & \\
& w_1 + w_2 = 1 & \\
& 1 \leq q \leq p \\
\text{where} \quad & k, q, n, x_i \in \mathbb Z_{> 0}, \\
\quad & f_i(x_i), g_i(x_i), \epsilon_1, \epsilon_2 \in \mathbb R_{> 0}, w_1, w_2 \in \mathbb R_{\geq 0}
\end{alignat*}

The vector, $z^* = (z_1^*,z_2^*)^T$, is the \textit{ideal} objective vector obtained by minimizing each of the objective functions individually subject to the constraints. If $d=1$, the sum of weighted deviations is minimized and the problem is similar to the weighting method \cite{Miettinen1999}.

\subsection{Informal Description of GBOPA} \label{app:gbopa}

\begin{figure}[!htbp]
	\centering
	\includegraphics[width=2.2in]{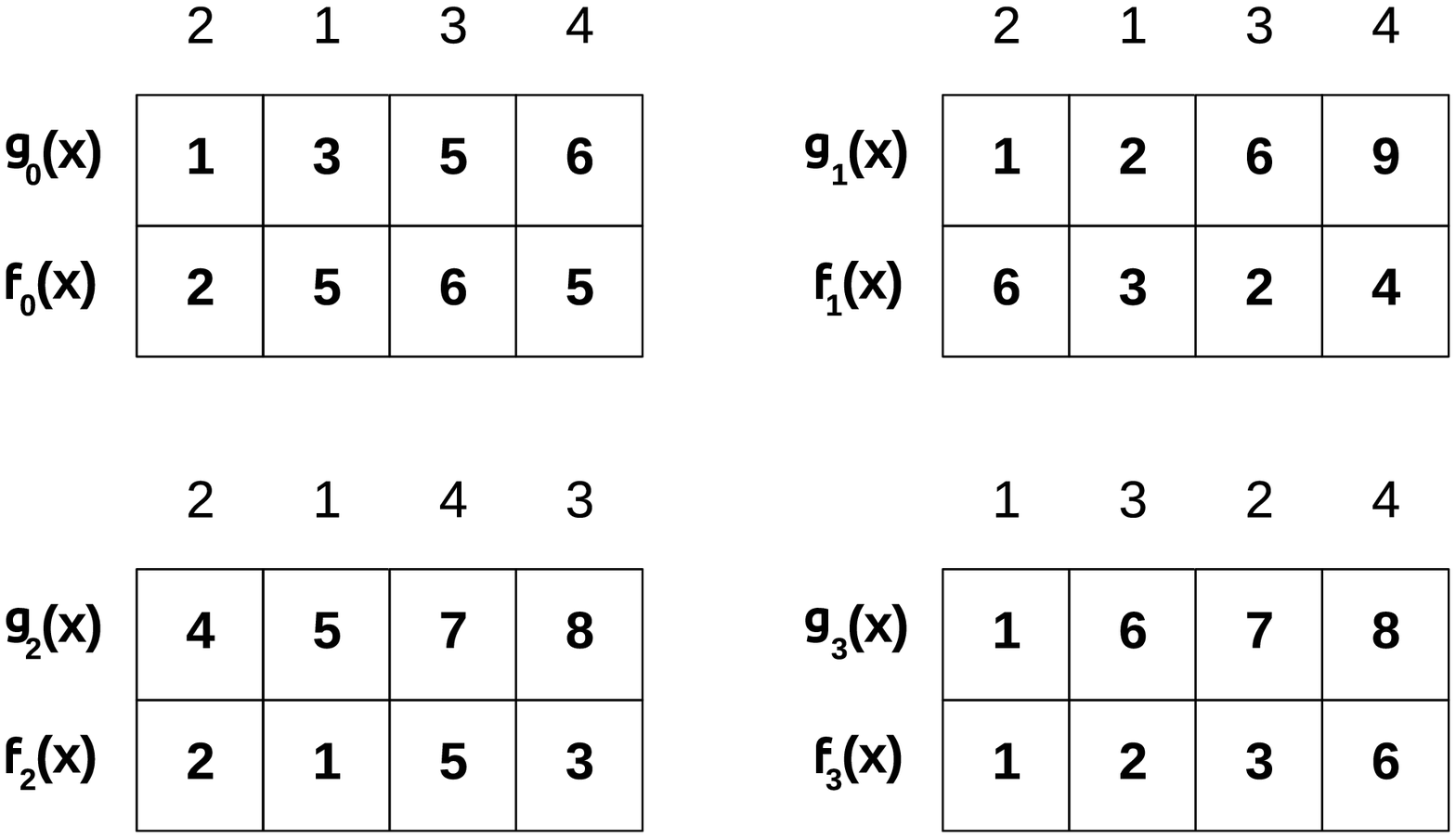}
	\caption{Functions in $F$ and $G$ given by sets of four points. The sets $f_i(x), g_i(x), i \in \{0,1,2,3\}$ are sorted in non-decreasing order of values for $g_i(x)$.}
	\label{fig:sample_functions}
\end{figure}

GBOPA is first described using a simple example where $k$ is equal to 4. Each function in the sets $F$ and $G$ is represented by a set of four points. The input $n$ is equal to 4. Figure \ref{fig:sample_functions} shows the functions. For instance, $g_0(1) = 3$ and $f_0(1) = 5$. The sets $f_i(x), g_i(x), i \in \{0,1,2,3\}$ are sorted in non-decreasing order of values for $g_i(x)$. The domain of $g_i(x)$ is represented by $D_i$.

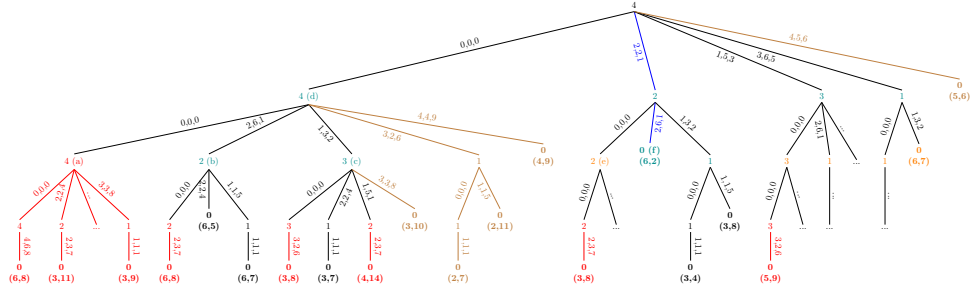
\begin{figure*}[!htbp]
	\centering
	\begin{tikzpicture} [scale=0.345,sloped]
	\tikzset{level 1/.style={level distance=3.5cm,sibling distance=0.5cm}}
	\tikzset{level 2/.style={level distance=2.5cm,sibling distance=0.5cm}}
	\tikzset{level 3/.style={level distance=2.5cm,sibling distance=0.45cm}}
	\tikzset{level 4/.style={level distance=2cm,sibling distance=0.1cm}}
		\Tree	[.4 \edge node[anchor=south] {0,0,0};
					[.\node[darkcyan] {4 (d)}; \edge node[anchor=south] {0,0,0};
						[.\node[red] {4 (a)}; \edge[red] node[anchor=south] {0,0,0};
							[.\node[red] {4}; \edge[red] node[anchor=south] {4,6,8};
								\node[red] [align=center]{\textbf{0}\\\textbf{(6,8)}};
							]
							\edge[red] node[anchor=south] {2,2,4};
							[.\node[red] {2}; \edge[red] node[anchor=south] {2,3,7};
								\node[red] [align=center]{\textbf{0}\\\textbf{(3,11)}};
							]
							\edge[red] node[anchor=south] {...};
							[ 
								.\node [red]{...};
							]
							\edge[red] node[anchor=south] {3,3,8};
							[.\node[red] {1}; \edge[red] node[anchor=south] {1,1,1};
								\node[red] [align=center]{\textbf{0}\\\textbf{(3,9)}};
							]
						]
						\edge node[anchor=south] {2,6,1};
						[.\node[darkcyan] {2 (b)};	\edge node[anchor=south] {0,0,0};
							[.\node[red] {2}; \edge[red] node[anchor=south] {2,3,7};
								\node [align=center,red]{\textbf{0}\\\textbf{(6,8)}};
							]
							\edge node[anchor=north] {2,2,4};
							[
								.\node [align=center]{\textbf{0}\\\textbf{(6,5)}};
							]
							\edge node[anchor=south] {1,1,5};
							[.1 \edge node[anchor=south] {1,1,1};
								\node [align=center]{\textbf{0}\\\textbf{(6,7)}};
							]
						]
						\edge node[anchor=north] {1,3,2};
						[.\node[darkcyan] {3 (c)}; \edge node[anchor=south] {0,0,0};
							[.\node[red] {3}; \edge[red] node[anchor=south] {3,2,6};
								\node [align=center,red]{\textbf{0}\\\textbf{(3,8)}};						
							]
							\edge node[anchor=north] {2,2,4};
							[.1 \edge node[anchor=south] {1,1,1};
								\node [align=center]{\textbf{0}\\\textbf{(3,7)}};						
							]
							\edge node[anchor=south] {1,5,1};
							[.\node[red] {2}; \edge[red] node[anchor=south] {2,3,7};
								\node [align=center,red]{\textbf{0}\\\textbf{(4,14)}};						
							]
							\edge[brown] node[anchor=south,brown] {3,3,8};
							[
								.\node[align=center,brown]{\textbf{0}\\\textbf{(3,10)}};							
							]					
						]
						\edge[brown] node[anchor=north] {3,2,6};
						[.\node[brown] {1}; \edge[brown] node[anchor=south,brown] {0,0,0};
							[.\node[brown] {1}; \edge[brown] node[anchor=south,brown] {1,1,1};
								\node[brown] [align=center,brown]{\textbf{0}\\\textbf{(2,7)}};						
							]
							\edge[brown] node[anchor=north,brown] {1,1,5};
							[
								.\node [align=center,brown]{\textbf{0}\\\textbf{(2,11)}};						
							]							
						]	
						\edge[brown] node[anchor=south,brown] {4,4,9};
						[
							.\node [align=center,brown]{\textbf{0}\\\textbf{(4,9)}};
						]
					]
					\edge[blue] node[anchor=north] {2,2,1};
					[.\node[darkcyan] {2};  \edge node[anchor=south] {0,0,0};
						[.\node[orange] {2 (e)}; \edge node[anchor=south] {0,0,0};
							[.\node[red]{2}; \edge[red] node[anchor=south] {2,3,7};
								\node [align=center,red]{\textbf{0}\\\textbf{(3,8)}};
							]
							\edge node[anchor=south] {...};
							[
								.\node {...};
							]
						]
						\edge[blue] node[anchor=north] {2,6,1};
						[
							.\node [align=center,color=darkcyan]{\textbf{0 (f)}\\\textbf{(6,2)}};
						]						
						\edge node[anchor=south] {1,3,2};
						[.\node[darkcyan] {1}; \edge node[anchor=south] {0,0,0};
							[.1 \edge node[anchor=south] {1,1,1};
								\node [align=center]{\textbf{0}\\\textbf{(3,4)}};
							]
							\edge node[anchor=south] {1,1,5};
							[
								.\node [align=center]{\textbf{0}\\\textbf{(3,8)}};
							]
						]
					]
					\edge node[anchor=north] {1,5,3};
					[.\node[darkcyan] {3}; \edge node[anchor=south] {0,0,0};
						[.\node[orange] {3}; \edge node[anchor=south] {0,0,0};
							[.\node[red] {3}; \edge[red] node[anchor=south] {3,2,6};
								\node [align=center,red]{\textbf{0}\\\textbf{(5,9)}};
							]
							\edge node[anchor=south] {...};
							[
								.\node {...};
							]
						]
						\edge node[anchor=north] {2,6,1};
						[.\node[orange] {1}; \edge node[anchor=south] {...};
							[
								.\node {...};
							]
						]
						\edge node[anchor=south] {...};
						[
							.\node {...};
						]
					]
					\edge node[anchor=north] {3,6,5};
					[.\node[darkcyan] {1}; \edge node[anchor=south] {0,0,0};
						[.\node[orange] {1}; \edge node[anchor=south] {...};
							[
								.\node {...};
							]				
						]
						\edge node[anchor=south] {1,3,2};
						[
							.\node [align=center,orange]{\textbf{0}\\\textbf{(6,7)}};				
						]
					]			
					\edge[brown] node[anchor=south,brown] {4,5,6};
					[
						.\node [align=center,brown]{\textbf{0}\\\textbf{(5,6)}};
					]
				]		
	\end{tikzpicture}
	\caption{The solution tree explored by the naive algorithm solving BOPGVECD to find all partitions of $n=4$ and its Pareto optimal solutions.}
	\label{fig:search_tree}
\end{figure*}

One straightforward approach to find the Pareto optimal solutions is to explore the full solution tree and find all possible partitions of $n$. Figure \ref{fig:search_tree} shows the tree. Only a partial tree is shown due to space constraints. 

The tree consists of 4 levels, $\{L_0,L_1,L_2,L_3\}$. All $x_i$s in $D_i$ are examined in level $L_i$. Each node in $L_i$, $i \in \{0,1,2,3\}$, is labelled by a size which is solved by GBOPA for the level $L_i$. Each edge connecting a node at level $L_i$ to its ancestor is labelled by a triple $(x_i, f_i(x_i), g_i(x_i))$ where $x_i \in D_i$.

The exploration process begins from the root to solve BOPGVECD for the input $n=4$ on the level $L_0$. Five sizes, four in $D_0 = \{2,1,3,4\}$ and zero, are examined one after another. Although there is no ordering assumption, we examine the data points in non-decreasing order of their values for $g_i(x)$. Examining the data points $\{0,2,1,3,4\}$ on $L_0$ expands the root into 5 children at $L_1$ representing the remaining size to be solved on level $L_1$. For instance, the edge $(2,2,1)$, highlighted in blue in Figure \ref{fig:search_tree}, indicates that $x_0$ is set to 2 at level $L_0$, where $f(2) = 2$ and $g(2) = 1$, and its child is labelled by 2 which equals the remaining size. 

In the same manner, each node in levels $\{L_1,L_2,L_3\}$ is expanded towards the leaves. Any leaf node labelled by 0 illustrates a solution, $X=\{x_0, x_1, x_2, x_3\}$ where the data points $x_0$, $x_1$, $x_2$ and $x_3$ are the first parts of the labels of the edges in the path from the root to the leaf. The objective values are shown below the solution node. For example, the blue path $\{(2,2,1),(2,6,1)\}$ in the tree highlights a solution $X = \{2, 2, 0, 0\}$, where its objective values are $MAX(4, F, X) = \max \{2,6,0,0\} = 6$ and $SUM(4, G, X) =  1 + 1 + 0 + 0 = 2$.

Due to lack of space, we have not shown the branches that do not provide any solution. In a non-solution branch, the summation of data points labelling the edges from the root to its leaf is greater than 4.

Each internal node that is not a leaf in the solution tree has either 5 children or one child that is a leaf. In the general case, there will be either $m+1$ children or one leaf child. There are two types of leaves. The \emph{solution} leaves are labelled by $0$ along with their objective values beneath. The \emph{no-solution} leaves are eliminated from the tree and, therefore, not shown. Each internal node at level $L_i$ is labelled by a positive number $w$. The node becomes a root of a solution tree for solving BOPGVECD with a problem size $w$ on level $L_i$ and is therefore constructed recursively. 

Once a solution is found, the algorithm updates the Pareto optimal set. In the end, the Pareto optimal set includes three members, (6,2,\{2,2,0,0\}), (3,4,\{2,1,0,1\}), (2,5,\{2,0,2,0\}), where each element, $(t, e,\{x_0,\cdots,x_3\})$, represents the Pareto optimal solution $X$ along with its objective values, $t = MIN(4,F,X)$ and $e = SUM(4,G,X)$.

The naive algorithm has exponential complexity. Therefore, we propose GBOPA, which is an efficient recursive algorithm to determine the Pareto optimal set of solutions. It has polynomial complexity. The algorithm shrinks the search space by utilizing two bounding criteria and also memorising the intermediate solutions to avoid exploring whole subtrees in the solution tree.

We will now describe how GBOPA efficiently solves the aforementioned example using the branch-and-bound technique. GBOPA starts building the search tree from its root in the depth-first manner. In a given level $L_i$, it examines all data points existing in $D_i$. GBOPA applies two bounding criteria to cut branches of the solution tree which lead to no solution. 

\begin{figure}[!htbp]
	\centering
	\includegraphics[width=2.5in]{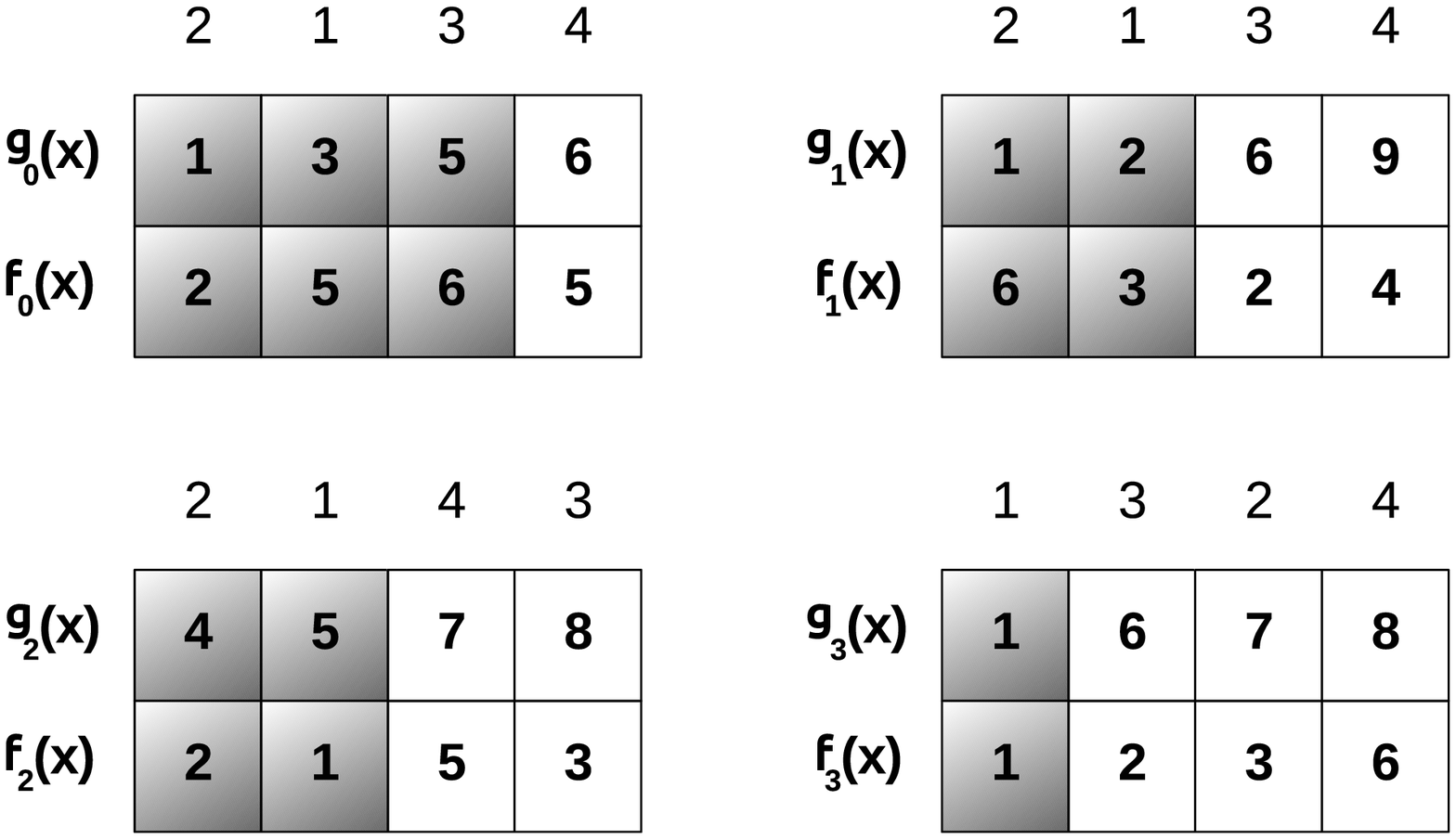}
	\caption{Removing data points from the sets by applying the sum threshold $\varepsilon$. The search space after the removal contains only the points in the shaded cells.}
	\label{fig:sample_functions_highlited}
\end{figure}

The first bounding criterion is called \emph{sum threshold}, represented by $\varepsilon$. The sum threshold $\varepsilon$ is initialised with $SUM(4,G,X_{min})$, where $X_{min}$ is the solution which optimizes the objective function, $MAX(4, F, X)$. We find this solution using the algorithm, HPOPTA \cite{khaleghzadeh2018novel} which solves the min-max optimisation problem. Applying sum threshold allows GBOPA to shrink search space by ignoring all data points $x$ with $g_i(x)$ greater than $\varepsilon$. In the example, the optimal solution, returned by HPOPTA, is $X_{min} = \{2,0,2,0\}$. Therefore, $\varepsilon$ in this example is set to 5 ($=SUM(4,G,X_{min})$). GBOPA, as shown in Figure \ref{fig:sample_functions_highlited}, ignores all the data points $x$ whose $g_i(x)$ is greater than 5. The search space after the removal contains only the points in the shaded cells. All branches eliminated from the solution tree by deploying sum threshold are highlighted in brown in Figure \ref{fig:search_tree}. There may exist more than one solution ($X_{min}$) minimizing the max objective but with different values for $\varepsilon = SUM(4,G,X_{min})$. The best solution has minimal $\varepsilon$. Nevertheless, using a non-optimal $\varepsilon$ does not restrain GBOPA from obtaining the set of Pareto optimal solutions.

To shrink the search space further, GBOPA assigns each level of the solution tree a size threshold $\sigma_i, i \in \{0,\dots,k-1\}$, which is considered as the second bounding criterion. For a given level $L_i$, its corresponding size threshold ($\sigma_i$) determines the maximum size which can be solved on level $L_i$ so that $\forall j \in \{i, \cdots, k-1\}$, $g_j(x_j)$ is not greater than $\varepsilon$ and is formulated as follows:
\begin{equation}	\label{eq:sigma}
\sigma_i =	\begin{cases}	
				supremum(x),~x \in D_i,~g_i(x) \le \varepsilon					& \quad i = k - 1 \\
				\sigma_{i+1} + supremum(x),~x \in D_i,~g_i(x) \le \varepsilon	& \quad 0 \le i < k - 1
			\end{cases}
\end{equation}
The vector of size thresholds is determined using the functions, $g_i(x), i \in \{0,1,2,3\}$, and $\varepsilon$. In this example, the maximum data points, $x$ whose corresponding $g_i(x)$ values are less than $\varepsilon=5$ will be 3, 2, 2, and 1 for the levels $L_0$,  $L_1$,  $L_2$ and $L_3$ respectively. Therefore, the size threshold vector $\sigma$ contains four elements, $\sigma = \{\sigma_0,\sigma_1,\sigma_2,\sigma_3\} = \{8,5,3,1\}$, giving the thresholds for the levels, $\{L_0,L_1,L_2,L_3\}$, respectively. Before expanding each node, GBOPA compares its problem size with its corresponding size threshold. If the problem size exceeds the size threshold, the node is not expanded since it results in a solution with $SUM(4,G,X)$ greater than $\varepsilon$. 

After determining the two bounding criteria $\varepsilon$ and $\sigma$, GBOPA explores the solution tree from its root in the left-to-right and depth-first order. It gradually constructs Pareto optimal solutions for any node using the sets of Pareto optimal solutions of its children. It first allocates zero data point to $x_0$ and $x_1$ (Figure \ref{fig:search_tree}). The remaining problem size at the level $L_2$ is 4 which is labelled by $4 (a)$ in the tree. Since the problem size 4 is greater than the corresponding size threshold $\sigma_2=3$, the node is not expanded further and is cut. This optimization is called operation \emph{Cut}. We highlight in red all sub-trees eliminated from the search space using the operation \emph{Cut}.

Returning to the tree exploration, GBOPA examines the next node $2 (b)$ at the level $L_2$. Expansion of this node results in two solutions solving problem size 2 on levels $L_2$ and $L_3$. GBOPA updates the Pareto optimal set for this node and saves the solution in an array called $PMem$. 

GBOPA memorizes Pareto optimal solutions for each node in levels $\{L_1, \allowbreak \cdots, \allowbreak L_{k-2}\}$. The information stored for a given solution $X' = \{x_i, x_{i+1}, \cdots, x_{k-1}\}$ found for node with a problem size $w$ at level $L_i$, $i \in \{1,\cdots,k-2\}$, is a quintuple $<t', e', x_i, x\#, key>$ where $t' = SUBMAX(i,k,F,X')$, $e' = SUBSUM(i,k,G,X')$, $x\#$ is the number of $x \in X'$ and greater than zero, and $key$ is set to $g_{i+1} (w-x_i)$ and points to a solution already saved in $PMem$ for a node at level $L_{i+1}$ and problem size $w-x_i$ . We call this Pareto optimal solution at level $L_{i+1}$ a \textit{partial solution} for the problem size $w$. This partial solution may not exist for some nodes, which in this case is represented by $\emptyset$. Since the values for $e'$ are unique in a Pareto optimal set, we use $key$ as a pointer to partial solutions. For each solution leaf in levels $\{L_1,\cdots,L_{k-2}\}$, like $0 (f)$ in Figure \ref{fig:search_tree}, GBOPA memorizes a solution $\{<0,0,0,0,\emptyset>\}$.

Thus, the information saved for the node $2 (b)$ is a Pareto optimal set including two members, $\{<2,4,2,1,\emptyset>,<1,6,1,2,\emptyset>\}$. We call this key operation, \emph{SavePareto}. Green nodes in the solution tree highlight ones whose Pareto optimal sets are saved. After $2 (b)$, the node $3 (c)$ is examined. The solution saved for this node is $\{<2,5,2,2,\emptyset>\}$. If there exist two solutions with the same values for $t'$ and $e'$, the solution with less value for $x\#$ will be selected as the Pareto optimal solution.

GBOPA then backtracks to node $4 (d)$ at level $L_1$ and builds its Pareto optimal set by merging Pareto optimal sets saved for its children, $2 (b)$ and  $3 (c)$. Consider the edge $(2,6,1)$ connecting the node $4 (d)$ to $2 (b)$. Merging this edge with the Pareto optimal set which has been already saved for $2 (b)$, $\{<2,4,2,1,\emptyset>,<1,6,1,2,\emptyset>\}$, results in one Pareto optimal solution for the node $4 (d)$ which is saved as the quintuple $<6,5,2,2,\textbf{4}>$. In this solution, the last element, $4$, which is highlighted in bold, points to its partial solution in the node $2 (b)$ at $L_2$, which is $\{<2,\textbf{4},2,1,\emptyset>\}$. Merging the edge $(1,3,2)$ with the Pareto optimal set for $3 (c)$, $\{<2,5,2,2,\emptyset>\}$, results in a new solution $\{<3,7,1,3,5>\}$. Therefore, the Pareto optimal set for the node $4 (d)$ is $\{<3,7,1,3,5>, <6,5,2,2,4>\}$, which is saved in the memory. 

After building and saving the Pareto optimal set of the node $4 (d)$, GBOPA visits the node $2 (e)$ at the level $L_2$. This node has already been explored, and therefore, its Pareto optimal set is retrieved from $PMem$. We call this key operation, \emph{ReadParetoMem}. The nodes whose solutions are retrieved from the memory are highlighted in orange.

After visiting the other remaining nodes, GBOPA backtracks to the root and builds the Pareto optimal solutions for the problem size $4$ at level $L_0$ using the Pareto optimal sets saved for its children. Then it terminates.

GBOPA thus deploys three key operations, which are a). \emph{Cut}, b). \emph{SavePareto}, and c). \emph{ReadParetoMem}, to efficiently explore solution trees and build the set of Pareto optimal solutions.

\subsection{Subroutines Employed in GBOPA for Solving BOPGVEC} \label{app:subroutines}

\subsubsection{Function \emph{SizeThresholdCalc}}

Algorithm \ref{alg_sizeTh_HDePOPTA} shows the pseudocode of the function \Call{SizeThresholdCalc}{}, which calculates the size threshold array, $\sigma$. First, it determines the size threshold of $L_{p-1}$ by finding the greatest problem size in the function $G_{p-1}$ whose functional value is less than or equal to $\varepsilon$ (Line \ref{alg_sizeTh_HDePOPTA_largestW1}). Then, it calculates $\sigma_i$, $i \in \{0,1,,\cdots,k-2\}$ where $\sigma_i$ is the summation of $\sigma_{i+1}$ with the greatest work-size in the function $G_{i}$ whose functional value is less than or equal to $\varepsilon$ (Lines \ref{alg_sizeTh_HDePOPTA_loop1}-\ref{alg_sizeTh_HDePOPTA_loop2}).

\begin{algorithm}[H]
	\scriptsize
	\caption{Algorithm Determining Size Thresholds} 		\label{alg_sizeTh_HDePOPTA}
	\begin{algorithmic}[1]	
		\Function{SizeThresholdCalc}{$k, G, \varepsilon, \sigma$}		
			\State $\sigma_{k-1}$ $\gets$ $\max_{j=0}^{m-1}\{x_{(k-1)~j}~|~g_{(k-1)~j} \le \varepsilon\}$		\label{alg_sizeTh_HDePOPTA_largestW1}
			\ForAll{$i = k-2$; $ i \ge 0$; $ i{-}{-}$}	\label{alg_sizeTh_HDePOPTA_loop1}
				\State $\sigma_i$ $\gets$ $\sigma_{i+1}+\max_{j=0}^{m-1}\{x_{ij}~|~g_{ij} \le \varepsilon\}$
			\EndFor		\label{alg_sizeTh_HDePOPTA_loop2}
			\State \Return $\sigma$
		\EndFunction		
	\end{algorithmic}
\end{algorithm}

\subsubsection{Function \emph{Cut}}

The function \Call{Cut}{} returns \emph{TRUE} if the input workload $n$ is greater than the input size threshold $\sigma$ (Algorithm \ref{alg_cut_code_HDePOPTA}).

\begin{algorithm}[H]
	\scriptsize
	\caption{Algorithm Cutting Search Tree using the Size Threshold} \label{alg_cut_code_HDePOPTA}
	\begin{algorithmic}[1]	
		\Function{Cut}{$n, \sigma$}
			\If{$n > \sigma$}
				\State \Return $TRUE$
			\EndIf
			\State\Return $FALSE$
		\EndFunction		
	\end{algorithmic}
\end{algorithm}

\subsubsection{Structure of matrix \emph{PMem}} \label{app:pmem}

We use \emph{PMem}, a two-dimensional array, to memorize Pareto-optimal solutions that have been found at levels $\{L_1,\cdots,L_{k-2}\}$ in the solution trees. Consider a given memory cell, $PMem[i][n]$, which saves a Pareto-optimal solution found for a given workload $n$ at Level $L_i$, $i \in \{1,\cdots,k-1\}$. The memory cell consists of a set where each element is a tuple, $<t, e, part, K\#, key>$, storing one Pareto-optimal solution. 

The field $t$ stores the value of the objective function $T$; $e$ stores the value of the objective function $E$; $part$ determines the problem size assigned to the node at level $L_i$; $K\#$ represents the number of problem sizes in the solution that are greater than 0. $key$ is the $E$ objective functional value of a saved Pareto-optimal solution, provided it exists, for a node at the level $L_{i+1}$ labelled by $n-part$ where this Pareto-optimal solution is the partial solution for the node $n$. Since functional values of objective $E$ are unique in Pareto fronts, we use this parameter for pointing to partial solutions. $key$ operates as a pointer to partial solutions.

Elements in Pareto fronts are sorted in increasing order of $E$ objective functional values. If there exists no Pareto-optimal solution for the workload $n$ on the level $i$, its corresponding memory cell, $PMem[i][n]$, will contain one tuple where its $e$ field is set to the constant value $\_NS$ (i.e., \emph{No\_Solution}).

\subsubsection{Function \emph{ReadParetoMem}}

Algorithm \ref{readmem_readM} illustrates the function \Call{ReadParetoMem}{}. Suppose we want to retrieve the saved solutions for a given workload $n$ on Level $L_{lvl}$. First, $PMem[lvl][n]$ is read, which contains the saved solutions for $n$ (Line \ref{readmem_readM}). If $PMem[lvl][n]$ is empty, which means that this node has not been visited yet, the function returns \emph{DUMMY} (Lines \ref{readmem_empty1}-\ref{readmem_empty2}). In this case, \emph{HDePOPTA\_Kernel} will continue with expanding this node.

Since solutions in memory cells are sorted in the increasing order of $E$ objective functional values, we consider the value of the first element in each set as the best solution. Based on the retrieved value for $e$, the following cases might happen:

\begin{itemize}
	\item \textbf{NOT\_SOLUTION}: This case occurs when $e$ is equal to $\_NS$ (there is no solution for $n$ on levels $\{L_{lvl},\cdots,L_{k-1}\}$) or the value of objective function $E$ of the saved solution is greater than $\varepsilon$ (Lines \ref{readmem_noSol_1} and \ref{readmem_noSol_2}).
	\item \textbf{SOLUTION}: This case occurs if the retrieved $e$ is less than or equal to $\varepsilon$ (Line \ref{readmem_Sol}).
\end{itemize}

\begin{algorithm}[H]
	\scriptsize
	\caption{Algorithm Retrieving Solution from Memory} \label{readmem_readM}
	\begin{algorithmic}[1]	
		\Function{ReadParetoMem}{$n, lvl, \varepsilon, PMem$}
			\State $pSet \gets PMem[lvl][n]$		\label{readmem_cSet}
			\If{$|pSet| = 0$}					\label{readmem_empty1}
				\State \Return $DUMMY$
			\EndIf								\label{readmem_empty2}
			\If{$pSet[0].e = \_NS \vee pSet[0].e > \varepsilon$}	\label{readmem_noSol_1}
				\State \Return $NOT\_SOLUTION$	
			\EndIf															\label{readmem_noSol_2}
			\State \Return $SOLUTION$						\label{readmem_Sol}
		\EndFunction		
	\end{algorithmic}
\end{algorithm}

\subsubsection{Function \emph{MakeParetoFinal}}

Algorithm \ref{alg_fin_HDePOPTA} illustrates the function \Call{MakeParetoFinal}{}, which finalizes the input memory cell, $pmem$. Each memory cell is finalized when its corresponding node and its children in the tree are completely explored. If a node is expanded for which there is no Pareto-optimal solution, the node is labelled as $\_NS$ by inserting a tuple with the constant value $\_NS$ in the field $e$ (Line \ref{finalize_insert}).

\begin{algorithm}[H]
	\scriptsize
	\caption{Algorithm Finalizing Memory Cells} \label{alg_fin_HDePOPTA}
	\begin{algorithmic}[1]	
		\Function{MakeParetoFinal}{$pmem$}
			\If{$|pmem| = 0$} 
				\State $pSet\gets(\_NS,0,0,0,0)$	\label{finalize_insert}	
			\EndIf
		\EndFunction		
	\end{algorithmic}
\end{algorithm}

\subsubsection{Function \emph{MergePartialParetoes}} \label{sec:merge}

For every non-leaf node, \emph{GBOPA\_Kernel} invokes the function \Call{MergePartialParetoes}{} to build its Pareto-optimal solutions, using the Pareto fronts of its children, which are named \emph{partial solutions} for the node. The function then stores the new solutions in $PMem$. If there exist two partitions with equal $T$ and $E$ objective functional values, \Call{MergePartialParetoes}{} selects the solution that has minimum number of problem sizes greater than 0. The input variable $lvl$ indicates a level in the tree, and $partsVec$ is a list including all problem sizes on level $L_{lvl}$ that result in a solution. The algorithm starts with initializing $pSet$, which points to a memory cell storing Pareto-optimal solutions for a workload $n$ on $L_{lvl}$ (Lines \ref{merge_pset1}-\ref{merge_pset2}). The set $\Psi$ will store final Pareto-optimal solutions for the root. The first $For$ loop iterates all problem sizes in $partsVec$ and builds new feasible solutions by merging the problem sizes in $partsVec$ with their corresponding partial Pareto-optimal solutions (Lines \ref{merge_for1_1}-\ref{merge_for1_2}). In each iteration, for a given problem size $x$, \Call{MergePartialParetoes}{} finds the partial Pareto-optimal solutions ($subPareto$) in $PMem$ ($1 \le lvl < k-2$) or builds it ($lvl = k-2$) (Lines \ref{merge_subPareto1}-\ref{merge_subPareto2}). The inner $For$ loop scans all Pareto-optimal solutions in $subPareto$. It merges the problem size $x$, given to $x_{lvl}$ at level $L_{lvl}$, with Pareto-solutions in $subPareto$ for levels $\{L_{lvl+1},\cdots,L_{p-1}\}$ in the tree (Lines \ref{merge_innerloop1}-\ref{merge_innerloop2}). For each merged solution, $pSet$ is examined to verify whether a Pareto-optimal solution exists in the set. If it is the case, $pSet$ is updated, and all solutions that are Pareto-optimal are eliminated. Therefore, for each newly merged solution $(t_x,e_x,x,K\#_x,key)$, the following situations may happen:

\begin{enumerate}
	\item $pSet$ is empty and the solution is inserted (Line \ref{merge_empty1}).
	\item There exists a solution in $pSet$ that its $e$ is equal to $e_x$. In this case the saved solution is updated if either $e_x$ is less than $e$ or $K\#_x$ is less than $K\#$ (Lines \ref{merge_equal1}-\ref{merge_equal2}).
	\item The $e_x$ of the merged Pareto-optimal solution is greater than ones in the $pSet$. The solution is inserted in case its $t_x$ is less than the last solution in $pSet$ (Lines \ref{merge_equal2}-\ref{merge_end2}).
	\item The $e_x$ of the merged Pareto-optimal solution is less than ones in the $pSet$. The solution is inserted in $pSet$ after eliminating all non-Pareto-optimal solutions (Lines \ref{merge_end2}-\ref{merge_begin2}).
	\item The $e_x$ of the merged solution is somewhere at the middle of $pSet$. In this case, the solution is inserted in $pSet$, and all non-Pareto-optimal solutions are removed (Lines \ref{merge_begin2}-\ref{merge_middle2}).
\end{enumerate}

It should be mentioned that the function \emph{lower\_bound} returns a pointer to the first element in the $pSet$ that its $e$ is greater than or equal to $e_x$.

The algorithm prevents further iteration in $pSet$ if the $T$ objective functional value of the partial Pareto-optimal solution that is evaluated last is less than or equal to the $T$ objective functional value of problem size $x$ on $P_{lvl}$. In fact, the further scanning of the $pSet$ will not lead to a Pareto-optimal solution. It is because all Pareto-optimal solutions are sorted in increasing order of $E$ objective functional values, which consequently implies that the $T$ objective functional values are decreasing in each set. Thus, all solutions built using the following elements in the $pSet$ will have the same $T$ functional value as $f_{lvl}(x)$ but with greater $E$ functional value.

Finally, the function \Call{BuildParetoSols}{} is called to obtain the workload distribution for each solution in $\Psi$ (Lines \ref{merge_buildsol1}-\ref{merge_buildsol2}).

\begin{algorithm}[H]
	\tiny
	\caption{Algorithm Merging Partial-Pareto Solutions} \label{hepopt_merge}
	\begin{algorithmic}[1]
		\Function{\textbf{MergePartialParetoes}}{$n, k, lvl, F, G,partsVec,PMem,\Psi$}
			\If{$lvl = 0$}							\label{merge_pset1}
				\State $pSet \gets PMem[0][0]$
			\Else
				\State $pSet \gets PMem[lvl][n]$
			\EndIf									\label{merge_pset2}
			\ForAll{$x \in partsVec$}				\label{merge_for1_1}
				\If{$lvl < k - 2$}					\label{merge_subPareto1}
					\State $subPareto \gets PMem[lvl + 1][n - x]$
				\Else
					\State $x' \gets n-x$
					\State $K\#_{x'} \gets (x = 0~?~0~:~1)$
					\State $t_{x'} \gets \Call{ReadFunc}{T_{k-1},x'}$
					\State $e_{x'} \gets \Call{ReadFunc}{E_{k-1},x'}$
					\State $subPareto \gets (e_{x'}, t_{x'},x',K\#_{x'},-)$
				\EndIf								\label{merge_subPareto2}
				\State $t_x \gets \Call{ReadFunc}{T_{lvl},x}$
				\State $K\#_x \gets (x = 0~?~0~:~1)$
				\ForAll{$tup \in subPareto$}				\label{merge_innerloop1}
					\State $e_x \gets tup.e + \Call{ReadFunc}{E_{lvl},x}$
					\State $t_x \gets Max (tup.t, t_x)$
					\State $K\#_x \gets K\#_x + K\#_{tup}$
					\State $key \gets tup.e$
					\If{$|pSet| = 0$}			
						\State $pSet \gets (e_x, t_x, x, K\#_x, key)$	\label{merge_empty1}
					\Else	
						\State $tup_l \gets pSet.lower\_bound(e_x)$
						\If{$tup_l \neq pSet.end() \wedge tup_l.e = e_x$}	\label{merge_equal1}
							\If{$tup_l.t > t_x$}
								\State $tup_l \gets (e_x, t_x, x, K\#_x, key)$
								\ForAll{$r \in pSet~|~r.e > e_x \wedge r.t \ge t_x$}
									\State $pSet \gets pSet - r$
								\EndFor
							\ElsIf{$tup_l.t = t_x \wedge K\#_x < tup_l.K\#$}
								\State $tup_l \gets (e_x, t_x, x, K\#_x, key)$
							\EndIf
						\ElsIf{$tup_l = pSet.end()$}						\label{merge_equal2}
							\State $tup_l \gets tup_l - 1$
							\If{$tup_l.t > t_x$}
								\State $pSet \cup (e_x, t_x, x, K\#_x, key)$
							\EndIf
						\ElsIf{$tup_l = pSet.begin()$}
							\If{$t_x \leq tup_l.t $}					\label{merge_end2}
								\ForAll{$r \in pSet~|~r.e > e_x \wedge r.t \ge t_x$}
									\State $pSet \gets pSet - r$
								\EndFor
							\EndIf
							\State $pSet \cup (e_x, t_x, x, K\#_x, key)$
						\Else												\label{merge_begin2}
							\State $tup_l \gets tup_l - 1$
							\If{$tup_l.t > t_x$}
								\State $pSet \cup (e_x, t_x, x, K\#_x, key)$
								\ForAll{$r \in pSet~|~r.e > e_x \wedge r.t \ge t_x$}
									\State $pSet \gets pSet - r$
								\EndFor
							\EndIf
						\EndIf												\label{merge_middle2}
						\If{$tup.t \leq t_x$}
							\State \textbf{break}							\label{merge_break}
						\EndIf
					\EndIf
				\EndFor		\label{merge_innerloop2}
			\EndFor		\label{merge_for1_2}
			\If{$c = 0$} \label{merge_buildsol1}
				\State \Call {BuildParetoSols}{$PMem, \Psi$}	\label{HDePOPTA_buildsols}
			\EndIf									\label{merge_buildsol2}
		\EndFunction
	\end{algorithmic}	
\end{algorithm}

\subsubsection{Function \emph{BuildParetoSols}}
As explained in the section \ref{sec:merge}, the set $\Psi$ holds final Pareto-optimal solutions. Each element in $\Psi$, which represents a Pareto-optimal solution, is a triple, $(e,t,X)$ where $e$ and $t$ are the $E$ and $T$ objective functional values, and $X=\{x_0, x_1, \cdots, x_{k-1}\}$ represents the partition. The function \Call{BuildParetoSols}{} determines the partition stored in $PMem[0][0]$. 

The algorithm \ref{alg_buildPareto} shows the pseudocode of \Call{BuildParetoSols}{}. The function reads the problem sizes in the solution partition from $PMem$. It uses the field $key$ in each saved solution to find the corresponding partial solution and eventually the problem size, $X_{i+1}$. Since $E$ objective functional values are unique in any set, there is only one tuple whose $E$ objective functional value is equal to $key$ in that set.

\begin{algorithm}[H]
	\scriptsize
	\caption{Algorithm Completing Workload Distribution for $\Psi$} \label{alg_buildPareto}
	\begin{algorithmic}[1]	
		\Function{BuildParetoSols}{$k, PMem, \Psi$}
			\ForAll{$tup \in PMem[0][0]$}
				\State $sumSize \gets tup.part$
				\State $X[0] \gets tup.part$
				\State $key_{cur} \gets tup.key$
				\ForAll{$i = 1$; $ i \le k-2$; $ i{+}{+}$}
					\State $tup_{sub} \gets \{t \in PMem[i][n - sumSize]~|~t.e = key_{cur}\}$
					\State $X[i] \gets tup_{sub}.part$
					\State $sumSize \gets sumSize + tup_{sub}.part$
					\State $key_{cur} \gets tup_{sub}.key$
				\EndFor
				\State $X[k-1] \gets n - sumSize$
				\State $\Psi \gets \Psi \cup (tup.e, tup.t, X)$
			\EndFor
		\EndFunction		
	\end{algorithmic}
\end{algorithm}

\subsubsection{Sample Solution Tree for Complexity Proof of GBOPA}

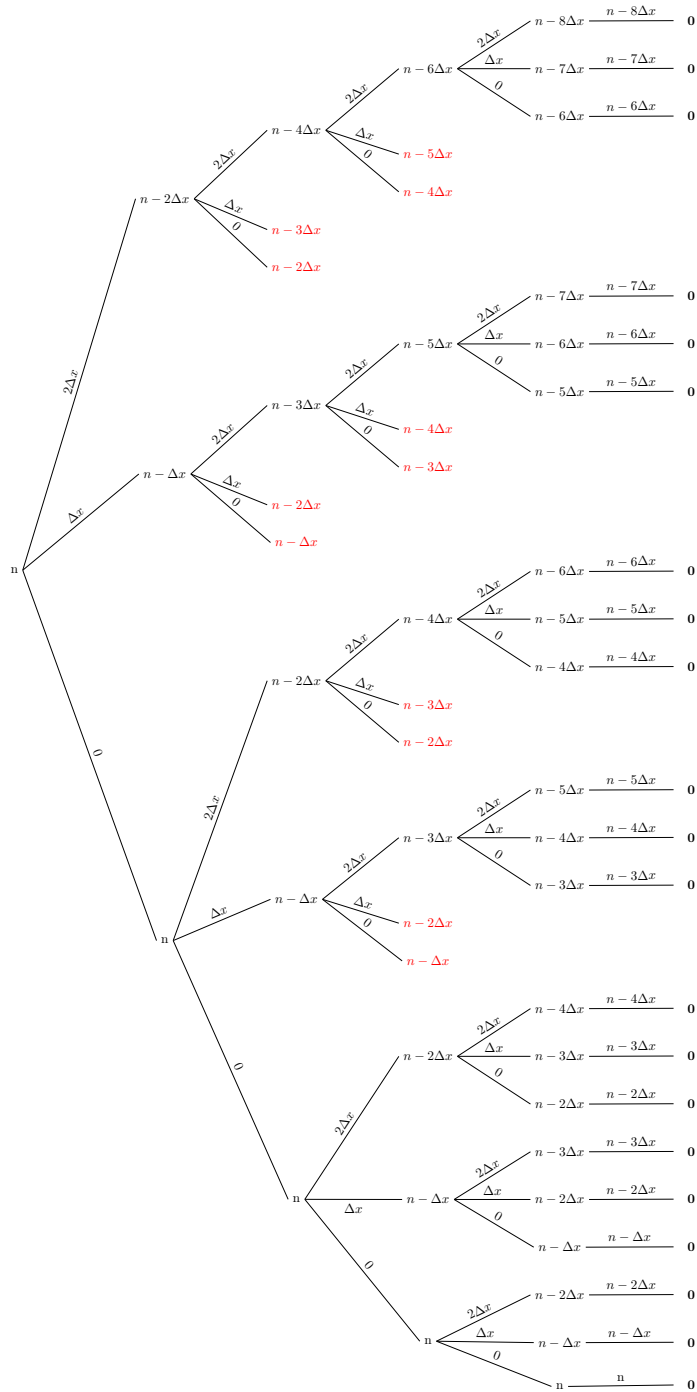
\begin{figure}[H]
	\centering
	\begin{tikzpicture} [scale=0.5,sloped,grow=right]
	\tikzset{level 1/.style={level distance=4cm,sibling distance=0cm}}
	\tikzset{level 2/.style={level distance=3.5cm,sibling distance=0.5cm}}
	\tikzset{level 3/.style={level distance=3.5cm,sibling distance=0.5cm}}
	\tikzset{level 4/.style={level distance=3.5cm,sibling distance=0.5cm}}
	\tikzset{level 5/.style={level distance=3.5cm,sibling distance=0cm}}
		\Tree	[.n \edge node[anchor=south] {0};
					[.n \edge node[anchor=south] {0};
						[.n \edge node[anchor=south] {0};
							[.n \edge node[anchor=south] {0};
								[.n \edge node[anchor=south] {n};
									\node [text width=0.7cm,align=center]{\textbf{0}};
								]
								\edge node[anchor=south] {$\Delta x$};
								[.$n-\Delta x$ \edge node[anchor=south] {$n-\Delta x$};
									\node [text width=0.7cm,align=center]{\textbf{0}};
								]
								\edge node[anchor=south] {$2\Delta x$};
								[.$n-2\Delta x$ \edge node[anchor=south] {$n-2\Delta x$};
									\node [text width=0.7cm,align=center]{\textbf{0}};
								]
							]
							\edge node[anchor=north] {$\Delta x$};
							[.$n-\Delta x$ \edge node[anchor=south] {0};
								[.$n-\Delta x$ \edge node[anchor=south] {$n-\Delta x$};
									\node [text width=0.7cm,align=center]{\textbf{0}};
								]
								\edge node[anchor=south] {$\Delta x$};
								[.$n-2\Delta x$ \edge node[anchor=south] {$n-2\Delta x$};
									\node [text width=0.7cm,align=center]{\textbf{0}};
								]
								\edge node[anchor=south] {$2\Delta x$};
								[.$n-3\Delta x$ \edge node[anchor=south] {$n-3\Delta x$};
									\node [text width=0.7cm,align=center]{\textbf{0}};
								]
							]
							\edge node[anchor=south] {$2\Delta x$};
							[.$n-2\Delta x$ \edge node[anchor=south] {0};
								[.$n-2\Delta x$ \edge node[anchor=south] {$n-2\Delta x$};
									\node [text width=0.7cm,align=center]{\textbf{0}};
								]
								\edge node[anchor=south] {$\Delta x$};
								[.$n-3\Delta x$ \edge node[anchor=south] {$n-3\Delta x$};
									\node [text width=0.7cm,align=center]{\textbf{0}};
								]
								\edge node[anchor=south] {$2\Delta x$};
								[.$n-4\Delta x$ \edge node[anchor=south] {$n-4\Delta x$};
									\node [text width=0.7cm,align=center]{\textbf{0}};
								]
							]
						]
						\edge node[anchor=south] {$\Delta x$};
						[.$n-\Delta x$ \edge node[anchor=south] {0};
							[.\node [red]{$n-\Delta x$};
							]
							\edge node[anchor=south] {$\Delta x$};
							[.\node [red]{$n-2\Delta x$};
							]
							\edge node[anchor=south] {$2\Delta x$};
							[.$n-3\Delta x$ \edge node[anchor=south] {0};
								[.$n-3\Delta x$ \edge node[anchor=south] {$n-3\Delta x$};
									\node [text width=0.7cm,align=center]{\textbf{0}};
								]
								\edge node[anchor=south] {$\Delta x$};
								[.$n-4\Delta x$ \edge node[anchor=south] {$n-4\Delta x$};
									\node [text width=0.7cm,align=center]{\textbf{0}};
								]
								\edge node[anchor=south] {$2\Delta x$};
								[.$n-5\Delta x$ \edge node[anchor=south] {$n-5\Delta x$};
									\node [text width=0.7cm,align=center]{\textbf{0}};
								]
							]
						]
						\edge node[anchor=south] {$2\Delta x$};
						[.$n-2\Delta x$ \edge node[anchor=south] {0};
							[.\node [red]{$n-2\Delta x$};
							]
							\edge node[anchor=south] {$\Delta x$};
							[.\node [red]{$n-3\Delta x$};
							]
							\edge node[anchor=south] {$2\Delta x$};
							[.$n-4\Delta x$ \edge node[anchor=south] {0};
								[.$n-4\Delta x$ \edge node[anchor=south] {$n-4\Delta x$};
									\node [text width=0.7cm,align=center]{\textbf{0}};
								]
								\edge node[anchor=south] {$\Delta x$};
								[.$n-5\Delta x$ \edge node[anchor=south] {$n-5\Delta x$};
									\node [text width=0.7cm,align=center]{\textbf{0}};
								]
								\edge node[anchor=south] {$2\Delta x$};
								[.$n-6\Delta x$ \edge node[anchor=south] {$n-6\Delta x$};
									\node [text width=0.7cm,align=center]{\textbf{0}};
								]
							]
						]
					]
					\edge node[anchor=south] {$\Delta x$};
					[.$n-\Delta x$ \edge node[anchor=south] {0};
						[.\node [red]{$n-\Delta x$};
						]
						\edge node[anchor=south] {$\Delta x$};
						[.\node [red]{$n-2\Delta x$};
						]
						\edge node[anchor=south] {$2\Delta x$};
						[.$n-3\Delta x$ \edge node[anchor=south] {0};
							[.\node [red]{$n-3\Delta x$};
							]
							\edge node[anchor=south] {$\Delta x$};
							[.\node [red]{$n-4\Delta x$};
							]
							\edge node[anchor=south] {$2\Delta x$};
							[.$n-5\Delta x$ \edge node[anchor=south] {0};
								[.$n-5\Delta x$ \edge node[anchor=south] {$n-5\Delta x$};
									\node [text width=0.7cm,align=center]{\textbf{0}};
								]
								\edge node[anchor=south] {$\Delta x$};
								[.$n-6\Delta x$ \edge node[anchor=south] {$n-6\Delta x$};
									\node [text width=0.7cm,align=center]{\textbf{0}};
								]
								\edge node[anchor=south] {$2\Delta x$};
								[.$n-7\Delta x$ \edge node[anchor=south] {$n-7\Delta x$};
									\node [text width=0.7cm,align=center]{\textbf{0}};
								]
							]
						]
					]
					\edge node[anchor=south] {$2\Delta x$};
					[.$n-2\Delta x$ \edge node[anchor=south] {0};
						[.\node [red]{$n-2\Delta x$};
						]
						\edge node[anchor=south] {$\Delta x$};
						[.\node [red]{$n-3\Delta x$};
						]
						\edge node[anchor=south] {$2\Delta x$};
						[.$n-4\Delta x$ \edge node[anchor=south] {0};
							[.\node [red]{$n-4\Delta x$};
							]
							\edge node[anchor=south] {$\Delta x$};
							[.\node [red]{$n-5\Delta x$};
							]
							\edge node[anchor=south] {$2\Delta x$};
							[.$n-6\Delta x$ \edge node[anchor=south] {0};
								[.$n-6\Delta x$ \edge node[anchor=south] {$n-6\Delta x$};
									\node [text width=0.7cm,align=center]{\textbf{0}};
								]
								\edge node[anchor=south] {$\Delta x$};
								[.$n-7\Delta x$ \edge node[anchor=south] {$n-7\Delta x$};
									\node [text width=0.7cm,align=center]{\textbf{0}};
								]
								\edge node[anchor=south] {$2\Delta x$};
								[.$n-8\Delta x$ \edge node[anchor=south] {$n-8\Delta x$};
									\node [text width=0.7cm,align=center]{\textbf{0}};
								]
							]
						]
					]	
				]		
	\end{tikzpicture}
    \caption{The GBOPA full solution tree for the bi-objective optimization problem BOPGVECD with $k = 5$ and $m = 2$. The memorization technique is only considered to reduce the full search space of solutions.}
	\label{fig:search_tree_time_com}
\end{figure}

Figure \ref{fig:search_tree_time_com} shows the sample full solution tree for the bi-objective optimization problem BOPGVECD with $k = 5$ and $m = 2$. Suppose $n$ is greater than $8\Delta x$, which is equal to the maximum possible size subtracted from $n$ in this example. In the figure, red nodes are ones that have been already expanded in the same level, and their solutions are retrieved from $PMem$. For the sake of simplicity, we only consider the two operations \emph{SavePareto} and \emph{ReadParetoMem}, and \emph{Cut} is ignored.

\subsection{Experimental Platform and Data-parallel Applications} \label{app:platform}

Our experimental platform consists of two heterogeneous nodes. The first node, HCLServer01, consists of an Intel Haswell multicore CPU (CPU\_1) involving 24 physical cores with 64 GB main memory, which hosts two accelerators, one Nvidia K40c GPU (GPU\_1) and one Intel Xeon Phi 3120P (Xeon Phi\_1) (specifications in Table \ref{table:hclserver12}). HCLServer02 contains an Intel Skylake multicore CPU (CPU\_2) consisting of 22 cores and 96 GB main memory. The multicore CPU is integrated with one Nvidia P100 GPU (GPU\_2) (specifications in Table \ref{table:hclserver12}). Each accelerator connects to a dedicated host core via a separate PCI-E link.

\begin{table}[!htbp]
	\caption{Specifications of the five heterogeneous processors.}
	\label{table:hclserver12}
	\centering
	\begin{tabular}{ |l|l| }
		\hline
		\multicolumn{2}{|c|}{\textbf{Intel Haswell E5-2670V3 (CPU\_1)}} \\ \hline
		No. of cores per socket & 12 \\ \hline
		Socket(s) & 2 \\ \hline
		CPU MHz & 1200.402 \\ \hline
		L1d cache, L1i cache  & 32 KB, 32 KB \\ \hline
		L2 cache, L3 cache & 256 KB, 30720 KB \\ \hline
		Total main memory & 64 GB DDR4 \\ \hline
		Memory bandwidth & 68 GB/sec \\ \hline
		\multicolumn{2}{|c|}{\textbf{NVIDIA K40c (GPU\_1)}} \\ \hline
		No. of processor cores & 2880 \\ \hline
		Total board memory & 12 GB GDDR5 \\ \hline
		L2 cache size & 1536 KB \\ \hline
		Memory bandwidth & 288 GB/sec \\ \hline
		\multicolumn{2}{|c|}{\textbf{Intel Xeon Phi 3120P (Xeon Phi\_1)}} \\ \hline
		No. of processor cores & 57 \\ \hline
		Total main memory & 6 GB GDDR5 \\ \hline
		Memory bandwidth & 240 GB/sec \\ \hline
		\multicolumn{2}{|c|}{\textbf{Intel Xeon Gold 6152 (CPU\_2)}} \\ \hline
		Socket(s) & 1 \\ \hline
		Cores per socket & 22 \\ \hline
		L1d cache, L1i cache  & 32 KB, 32 KB \\ \hline
		L2 cache, L3 cache & 256 KB, 30976 KB \\ \hline        
		Main memory &  96 GB \\ \hline
		\multicolumn{2}{|c|}{\textbf{NVIDIA P100 PCIe (GPU\_2)}} \\ \hline
		No. of processor cores & 3584 \\ \hline
		Total board memory & 12 GB CoWoS HBM2 \\ \hline
		Memory bandwidth & 549 GB/sec \\ \hline
	\end{tabular}
\end{table}

A data-parallel application executing on this heterogeneous hybrid platform, consists of several kernels (generally speaking, multithreaded), running in parallel on different computing devices of the platform. The proposed algorithm for solving the bi-objective optimisation problem for performance and energy requires individual performance and energy profiles of all the kernels. Due to tight integration and severe resource contention in heterogeneous hybrid platforms, the load of one computational kernel in a given hybrid application may significantly impact others' performance to the extent of preventing the ability to model the performance and energy consumption of each kernel in hybrid applications individually \cite{Zhong2015}. To address this issue, we restrict our study in this work to configurations of hybrid applications, where individual kernels are coupled loosely enough to allow us to build their performance and energy profiles with the accuracy sufficient for successful application of the proposed algorithms. To achieve this objective, we only consider configurations where no more than one CPU kernel or accelerator kernel runs on the corresponding device. To apply our optimization algorithms, each group of cores executing an individual kernel of the application is modelled as an abstract processor \cite{Zhong2015}, so that the executing platform is represented as a set of heterogeneous abstract processors. We make sure that the sharing of system resources is maximized within groups of computational cores representing the abstract processors and minimized between the groups. This way, the contention and mutual dependence between abstract processors are minimized.

We thus model HCLServer01 by three abstract processors, CPU\_1, GPU\_1, and PHI\_1. CPU\_1 represents 22 (out of total 24) CPU cores. GPU\_1 involves the Nvidia K40c GPU and a host CPU core connected to this GPU via a dedicated PCI-E link. PHI\_1 is made up of one Intel Xeon Phi 3120P and its host CPU core connected via a dedicated PCI-E link. In the same manner, HCLServer02 is modelled by two abstract processors, CPU\_2 and GPU\_2. Since there should be a one-to-one mapping between the abstract processors and computational kernels, any hybrid application executing on the servers in parallel should consist of five kernels, one kernel per computational device. Because the abstract processors contain CPU cores that share some resources such as main memory and QPI, they cannot be considered entirely independent. Therefore, the performance of these loosely-coupled abstract processors must be measured simultaneously, thereby taking into account the influence of resource contention \cite{Zhong2015}. 

The execution time of any computational kernel can be measured accurately using high precision processor clocks and used to model the performance of a parallel application and build its speed functions. There is, however, no such effective equivalent for measuring energy consumption. System-level physical measurements using power meters are accurate, but they do not provide a fine-grained decomposition of the energy consumption during the application run in a hybrid platform. Fahad et al. \cite{Fahad2020} propose a methodology to determine this decomposition, which employs only system-level power measurements using power meters. The methodology allows us to build discrete dynamic energy functions of abstract processors with sufficient accuracy for applying the proposed optimization algorithms in our use cases.

The matrix multiplication application, DGEMM, computes $C = \alpha \times A \times B + \beta \times C$, where $A$, $B$, and $C$ are matrices of size $m \times n$, $n \times n$, and $m \times n$, and $\alpha$ and $\beta$ are floating-point constants. The application uses Intel MKL DGEMM for CPUs, ZZGEMMOOC out-of-card package \cite{khaleghzadeh2018out} for Nvidia GPUs, and XeonPhiOOC out-of-card package \cite{khaleghzadeh2018out} for Intel Xeon Phis. ZZGEMMOOC and XeonPhiOOC packages reuse CUBLAS and MKL BLAS for in-card DGEMM calls. The out-of-card packages allow the GPUs and Xeon Phis to execute computations of arbitrary size. The Intel MKL and CUDA versions used on HCLServer01 are 2017.0.2 and 7.5, and on HCLServer02 are 2017.0.2 and 9.2.148. Workload sizes range from $64\times10112$ to $28800 \times10112$ with a step size of $64$ for the first dimension $m$. The speed of execution of a given problem size $m \times n$ is calculated as $(2 \times m \times n^2)/t$ where $t$ is the execution time.

The 2D fast Fourier transform application computes 2D-DFT of a complex signal matrix of size $m \times n$. It employs Intel MKL FFT routines for CPUs and Xeon Phis, and CUFFT routines for Nvidia GPUs. All computations are in-card. Workloads range from $1024\times51200$ to $10000\times51200$ with the step size of 16 for $m$. The experimental data set does not include problem sizes that cannot be factored into primes less than or equal to 127. For these problem sizes, CUFFT for GPU gives failures. The speed of execution of a 2D-DFT of size $m \times n$ is calculated as $(2.5 \times m \times n \times \log_2 (m \times n))/t$ where $t$ is the execution time.

\bibliographystyle{siamplain}
\bibliography{references}

\begin{thebibliography}{10}

\bibitem{aba2017approximation}
{\sc M.~A. Aba, L.~Zaourar, and A.~Munier}, {\em Approximation algorithm for
  scheduling a chain of tasks on heterogeneous systems}, in European Conference
  on Parallel Processing, Springer, 2017, pp.~353--365.

\bibitem{Berman1990}
{\sc O.~Berman, D.~Einav, and G.~Handler}, {\em The constrained bottleneck
  problem in networks}, Oper. Res., 38 (1990), p.~178–181.

\bibitem{Bornstein2012}
{\sc C.~T. Bornstein, N.~Maculan, M.~Pascoal, and L.~L. Pinto}, {\em
  Multiobjective combinatorial optimization problems with a cost and several
  bottleneck objective functions: An algorithm with reoptimization}, Comput.
  Oper. Res., 39 (2012), p.~1969–1976.

\bibitem{chakrabarti2017pareto}
{\sc A.~Chakrabarti, S.~Parthasarathy, and C.~Stewart}, {\em A pareto framework
  for data analytics on heterogeneous systems: Implications for green energy
  usage and performance}, in Parallel Processing (ICPP), 2017 46th
  International Conference on, IEEE, 2017, pp.~533--542.

\bibitem{Leizer2009}
{\sc L.~de~Lima~Pinto, C.~T. Bornstein, and N.~Maculan}, {\em The tricriterion
  shortest path problem with at least two bottleneck objective functions},
  European Journal of Operational Research, 198 (2009), pp.~387 -- 391.

\bibitem{Pinto2019}
{\sc L.~de~Lima~Pinto, K.~C.~C. Fernandes, K.~V. Cardoso, and N.~Maculan}, {\em
  An exact and polynomial approach for a bi-objective integer programming
  problem regarding network flow routing}, Comput. Oper. Res., 106 (2019),
  pp.~28--35.

\bibitem{Durillo2014}
{\sc J.~J. Durillo, V.~Nae, and R.~Prodan}, {\em Multi-objective
  energy-efficient workflow scheduling using list-based heuristics}, Future
  Generation Computer Systems, 36 (2014), pp.~221 -- 236.

\bibitem{fahad2019comparative}
{\sc M.~Fahad, A.~Shahid, R.~R. Manumachu, and A.~Lastovetsky}, {\em A
  comparative study of methods for measurement of energy of computing},
  Energies, 12 (2019), p.~2204.

\bibitem{Fahad2020}
{\sc M.~{Fahad}, A.~{Shahid}, R.~R. {Manumachu}, and A.~{Lastovetsky}}, {\em
  Accurate energy modelling of hybrid parallel applications on modern
  heterogeneous computing platforms using system-level measurements}, IEEE
  Access, 8 (2020), pp.~93793--93829.

\bibitem{gholkar2016power}
{\sc N.~Gholkar, F.~Mueller, and B.~Rountree}, {\em Power tuning {HPC} jobs on
  power-constrained systems}, in Proceedings of the 2016 International
  Conference on Parallel Architectures and Compilation, ACM, 2016,
  pp.~179--191.

\bibitem{Hansen1980}
{\sc P.~Hansen}, {\em Bicriterion path problems}, in Multiple Criteria Decision
  Making Theory and Application, G.~Fandel and T.~Gal, eds., Berlin,
  Heidelberg, 1980, Springer Berlin Heidelberg, pp.~109--127.

\bibitem{heydari2018dynamic}
{\sc J.~Heydari, A.~Sabbaghnia, and J.~Razmi}, {\em A dynamic bi-objective
  model for after disaster blood supply chain network design; a robust
  possibilistic programming approach}, Journal of Industrial and Systems
  Engineering, 11 (2018), pp.~16--28.

\bibitem{karimi2010bi}
{\sc N.~Karimi, M.~Zandieh, and H.~Karamooz}, {\em Bi-objective group
  scheduling in hybrid flexible flowshop: a multi-phase approach}, Expert
  Systems with Applications, 37 (2010), pp.~4024--4032.

\bibitem{Kessaci2013}
{\sc Y.~Kessaci, N.~Melab, and E.-G. Talbi}, {\em A pareto-based metaheuristic
  for scheduling {HPC} applications on a geographically distributed cloud
  federation}, Cluster Computing, 16 (2013), pp.~451--468.

\bibitem{HamidTPDS2020}
{\sc H.~{Khaleghzadeh}, M.~{Fahad}, A.~{Shahid}, R.~R. {Manumachu}, and
  A.~{Lastovetsky}}, {\em Bi-objective optimization of data-parallel
  applications on heterogeneous {HPC} platforms for performance and energy
  through workload distribution}, IEEE Transactions on Parallel and Distributed
  Systems, 32 (2021), pp.~543--560.

\bibitem{khaleghzadeh2018novel}
{\sc H.~Khaleghzadeh, R.~R. Manumachu, and A.~Lastovetsky}, {\em A novel
  data-partitioning algorithm for performance optimization of data-parallel
  applications on heterogeneous hpc platforms}, IEEE Transactions on Parallel
  and Distributed Systems, 29 (2018), pp.~2176--2190.

\bibitem{khaleghzadeh2018out}
{\sc H.~Khaleghzadeh, Z.~Zhong, R.~Reddy, and A.~Lastovetsky}, {\em Out-of-core
  implementation for accelerator kernels on heterogeneous clouds}, The Journal
  of Supercomputing, 74 (2018), pp.~551--568.

\bibitem{Kolodziej2015}
{\sc J.~Ko{\l}odziej, S.~U. Khan, L.~Wang, and A.~Y. Zomaya}, {\em Energy
  efficient genetic-based schedulers in computational grids}, Concurr. Comput.
  : Pract. Exper., 27 (2015), pp.~809--829.

\bibitem{Lang2014}
{\sc J.~Lang and G.~Rünger}, {\em An execution time and energy model for an
  energy-aware execution of a conjugate gradient method with {CPU/GPU}
  collaboration}, Journal of Parallel and Distributed Computing, 74 (2014),
  pp.~2884 -- 2897.

\bibitem{Lastovetsky2004}
{\sc A.~{Lastovetsky} and R.~{Reddy}}, {\em Data partitioning with a realistic
  performance model of networks of heterogeneous computers}, in 18th
  International Parallel and Distributed Processing Symposium, 2004.
  Proceedings., 2004, pp.~104--.

\bibitem{Lastovetsky2007}
{\sc A.~Lastovetsky and R.~Reddy}, {\em Data partitioning with a functional
  performance model of heterogeneous processors}, International Journal of High
  Performance Computing Applications, 21 (2007), pp.~76--90.

\bibitem{LastovetskyReddy2017}
{\sc A.~Lastovetsky and R.~Reddy}, {\em New model-based methods and algorithms
  for performance and energy optimization of data parallel applications on
  homogeneous multicore clusters}, IEEE Transactions on Parallel and
  Distributed Systems, 28 (2017), pp.~1119--1133.

\bibitem{manumachu2018bi}
{\sc R.~R. Manumachu and A.~Lastovetsky}, {\em Bi-objective optimization of
  data-parallel applications on homogeneous multicore clusters for performance
  and energy}, IEEE Transactions on Computers, 67 (2018), pp.~160--177.

\bibitem{melamed1996computational}
{\sc I.~Melamed and I.~K. Sigal}, {\em A computational investigation of linear
  parametricization of criteria in multicriteria discrete programming},
  Computational mathematics and mathematical physics, 10 (1996),
  pp.~1341--1343.

\bibitem{melamed1998numerical}
{\sc I.~Melamed and I.~K. Sigal}, {\em Numerical analysis of tricriteria tree
  and assignment problems}, Computational mathematics and mathematical physics,
  38 (1998), pp.~1707--1714.

\bibitem{Miettinen1999}
{\sc K.~Miettinen}, {\em Nonlinear multiobjective optimization}, Kluwer, 1999.

\bibitem{minoux1989solving}
{\sc M.~Minoux}, {\em Solving combinatorial problems with combined
  min-max-min-sum objective and applications}, Mathematical Programming, 45
  (1989), pp.~361--372.

\bibitem{Leizer2010}
{\sc L.~L. Pinto and M.~M. Pascoal}, {\em On algorithms for the tricriteria
  shortest path problem with two bottleneck objective functions}, Computers \&
  Operations Research, 37 (2010), pp.~1774 -- 1779.

\bibitem{Punnen1994}
{\sc A.~P. Punnen}, {\em On combined minmax-minsum optimization}, Computers \&
  Operations Research, 21 (1994), pp.~707 -- 716.

\bibitem{Punnen1996}
{\sc A.~P. Punnen and K.~Nair}, {\em An o(m log n) algorithm for the max + sum
  spanning tree problem}, European Journal of Operational Research, 89 (1996),
  pp.~423 -- 426.

\bibitem{manumachu2018bicpe}
{\sc R.~Reddy~Manumachu and A.~L. Lastovetsky}, {\em Design of self-adaptable
  data parallel applications on multicore clusters automatically optimized for
  performance and energy through load distribution}, Concurrency and
  Computation: Practice and Experience, 31 (2019), p.~e4958.

\bibitem{Rountree2017}
{\sc B.~{Rountree}, D.~K. {Lowenthal}, S.~{Funk}, V.~W. {Freeh}, B.~R. {de
  Supinski}, and M.~{Schulz}}, {\em Bounding energy consumption in large-scale
  {MPI} programs}, in SC '07: Proceedings of the 2007 ACM/IEEE Conference on
  Supercomputing, Nov 2007, pp.~1--9.

\bibitem{Ruzika2009}
{\sc S.~Ruzika and H.~W. Hamacher}, {\em A Survey on Multiple Objective Minimum
  Spanning Tree Problems}, Springer Berlin Heidelberg, Berlin, Heidelberg,
  2009, pp.~104--116.

\bibitem{salamati2018trade}
{\sc H.~Salamati-Hormozi, Z.-H. Zhang, O.~Zarei, and R.~Ramezanian}, {\em
  Trade-off between the costs and the fairness for a collaborative production
  planning problem in make-to-order manufacturing}, Computers \& Industrial
  Engineering, 126 (2018), pp.~421--434.

\bibitem{sergienko1987finding}
{\sc I.~Sergienko and V.~Perepelitsa}, {\em Finding the set of alternatives in
  discrete multicriterion problems}, Cybernetics, 23 (1987), pp.~673--683.

\bibitem{Sun2011}
{\sc Y.~Sun, C.~Zhang, L.~Gao, and X.~Wang}, {\em Multi-objective optimization
  algorithms for flow shop scheduling problem: a review and prospects}, The
  International Journal of Advanced Manufacturing Technology, 55 (2011),
  pp.~723--739.

\bibitem{Talbi2009}
{\sc E.-G. Talbi}, {\em Metaheuristics: from design to implementation},
  vol.~74, John Wiley \& Sons, 2009.

\bibitem{tarplee2016energy}
{\sc K.~M. Tarplee, R.~Friese, A.~A. Maciejewski, H.~J. Siegel, and E.~K.
  Chong}, {\em Energy and makespan tradeoffs in heterogeneous computing systems
  using efficient linear programming techniques}, IEEE Transactions on Parallel
  and Distributed Systems, 27 (2016), pp.~1633--1646.

\bibitem{Torkashvand2017}
{\sc M.~Torkashvand, B.~Naderi, and S.~Hosseini}, {\em Modelling and scheduling
  multi-objective flow shop problems with interfering jobs}, Appl. Soft
  Comput., 54 (2017), p.~221–228.

\bibitem{Yu2015}
{\sc L.~Yu, Z.~Zhou, S.~Wallace, M.~E. Papka, and Z.~Lan}, {\em Quantitative
  modeling of power performance tradeoffs on extreme scale systems}, Journal of
  Parallel and Distributed Computing, 84 (2015), pp.~1 -- 14.

\bibitem{Zhong2015}
{\sc Z.~Zhong, V.~Rychkov, and A.~Lastovetsky}, {\em Data partitioning on
  multicore and multi-{GPU} platforms using functional performance models},
  Computers, IEEE Transactions on, 64 (2015), pp.~2506--2518.

\end{thebibliography}
\end{document}